\documentclass{amsart}
\usepackage{times,amssymb,amscd,verbatim,color,graphics}
\usepackage[vcentermath,stdtext,enableskew]{youngtab}
\usepackage[all]{xy}

\numberwithin{equation}{section}

\newtheorem{prop}{Proposition}
\newtheorem{theorem}[prop]{Theorem}
\newtheorem{corollary}[prop]{Corollary}
\newtheorem{lemma}[prop]{Lemma}

\theoremstyle{definition}
\newtheorem{definition}[prop]{Definition}
\newtheorem{example}[prop]{Example}
\newtheorem{remark}[prop]{Remark}

\newtheorem{hypothesis}[prop]{Hypothesis}

\numberwithin{prop}{section}

\newcommand{\PBase}{LM}
\newcommand{\PBBaseA}{BC1}
\newcommand{\PBBaseB}{BC2}

\newcommand{\Conf}{\mathrm{C}}

\newcommand{\Confb}{\overline{\Conf}}

\newcommand{\J}{\overline{I}}
\newcommand{\RCb}{\overline{\RC}}
\newcommand{\D}{\mathcal{D}}
\newcommand{\Dom}{\mathrm{Dom}}

\newcommand{\HH}{\mathcal{H}}

\newcommand{\Img}{\mathrm{Img}}
\newcommand{\id}{\mathrm{id}}

\newcommand{\Plift}{l}
\newcommand{\Rlift}{\overline{l}}
\newcommand{\la}{\lambda}
\newcommand{\La}{\Lambda}
\newcommand{\lb}{\mathrm{lb}}

\newcommand{\Plb}{\mathrm{lb}}
\newcommand{\Psit}{\widetilde{\Psi}}
\newcommand{\Rlb}{\overline{\mathrm{lb}}}
\newcommand{\lh}{\mathrm{lh}}
\newcommand{\Plh}{\mathrm{lh}}
\newcommand{\Rlh}{\overline{\mathrm{lh}}}
\newcommand{\lm}{\la^-}

\newcommand{\ls}{\mathrm{ls}}
\newcommand{\Pls}{\mathrm{ls}}
\newcommand{\Rls}{\overline{\mathrm{ls}}}

\newcommand{\Path}{\mathcal{P}}
\newcommand{\Pathb}{\overline{\Path}}

\newcommand{\PRmap}{\mathrm{R}}

\newcommand{\RC}{\mathrm{RC}}
\newcommand{\Ppr}{\mathrm{pr}}
\newcommand{\Rpr}{\overline{\mathrm{pr}}}
\newcommand{\Pslide}{\rho}
\newcommand{\Rslide}{\overline{\rho}}
\newcommand{\Aslide}{\overline{\overline{\rho}}}
\newcommand{\rc}{\mathrm{rc}}
\newcommand{\rk}{\mathrm{rk}}
\newcommand{\rs}{\mathrm{rs}}
\newcommand{\Prs}{\mathrm{rs}}
\newcommand{\Rrs}{\overline{\mathrm{rs}}}

\newcommand{\wt}{\mathrm{wt}}
\newcommand{\Z}{\mathbb{Z}}

\newcommand{\row}{\mathrm{row}}
\newcommand{\col}{\mathrm{col}}

\newcommand{\Loneone}{S_{1,1}}

\newcommand{\Lonec}{S_{1,c}}

\newcommand{\Lrone}{S_{r,1}}
\newcommand{\Lrtwo}{S_{r,2}}

\newcommand{\Lrc}{S_{r,c}}
\newcommand{\Toneone}{T_{1,1}}

\newcommand{\Ttone}{T_{t,1}}

\newcommand{\two}{\textcolor{red}{2}}
\newcommand{\three}{\textcolor{red}{3}}
\newcommand{\four}{\textcolor{red}{4}}
\newcommand{\five}{\textcolor{red}{5}}
\newcommand{\six}{\textcolor{red}{6}}
\newcommand{\seven}{\textcolor{red}{7}}
\newcommand{\eight}{\textcolor{red}{8}}

\newcommand{\Tonek}{T_{1,s}}
\newcommand{\Trone}{T_{r,1}}

\newcommand{\Trk}{T_{r,s}}
\newcommand{\bONE}{b_1}
\newcommand{\bJ}{b_j}
\newcommand{\bR}{b_r}

\newcommand{\sONE}{s_1}
\newcommand{\sTWO}{s_2}
\newcommand{\sI}{s_k}

\newcommand{\sT}{s_t}

\begin{document}
\title{Promotion operator on rigged configurations of type $A$}

\author[A.~Schilling]{Anne Schilling}
\address{Department of Mathematics, University of California, One Shields
Avenue, Davis, CA 95616-8633, U.S.A.}
\email{anne@math.ucdavis.edu}
\urladdr{http://www.math.ucdavis.edu/\~{}anne}

\author[Q.~Wang]{Qiang Wang}
\email{xqwang@math.ucdavis.edu}
\urladdr{http://www.math.ucdavis.edu/\~{}xqwang}

\thanks{\textit{Date:} August 2009}
\thanks{Partially supported by NSF grants DMS--0501101, DMS--0652641, and DMS--0652652.}

\begin{abstract}
In~\cite{S:2006}, the analogue of the promotion operator on crystals of type $A$
under a generalization of the bijection of Kerov, Kirillov and Reshetikhin between crystals 
(or Littlewood--Richardson tableaux) and rigged configurations was proposed.
In this paper, we give a proof of this conjecture. This shows in particular that the bijection
between tensor products of type $A_n^{(1)}$ crystals and (unrestricted) rigged configurations is an 
affine crystal isomorphism.
\end{abstract}

\maketitle

%%%%%%%%%%%%%%%%%%%%%%%%%%%%%%%%%%%%%%%%%%%%%%%%%%%%%%%%%%%%%%%%%%%%%%%%%%%%%%%%%%%%%%%%%%%%%%%%%%%%%%%%%%%%%%%%%%%%%%%%%
%Section 1
%%%%%%%%%%%%%%%%%%%%%%%%%%%%%%%%%%%%%%%%%%%%%%%%%%%%%%%%%%%%%%%%%%%%%%%%%%%%%%%%%%%%%%%%%%%%%%%%%%%%%%%%%%%%%%%%%%%%%%%%%
\section{Introduction}
Rigged configurations appear in the Bethe Ansatz study of exactly solvable lattice models
as combinatorial objects to index the solutions of the Bethe equations~\cite{KKR:1986,KR:1988}.
Based on work by Kerov, Kirillov and Reshetikhin~\cite{KKR:1986, KR:1988}, it was shown
in~\cite{KSS:2002} that there is a statistic preserving bijection $\Phi$ between 
Littlewood-Richardson tableaux and rigged configurations.
The description of the bijection $\Phi$ is based on a quite technical recursive algorithm.

Littlewood-Richardson tableaux can be viewed as highest weight crystal elements in  a tensor 
product of Kirillov--Reshetikhin (KR) crystals of type $A_n^{(1)}$. KR crystals are
affine finite-dimensional crystals corresponding to affine Kac--Moody algebras, in the setting
of~\cite{KSS:2002} of type $A_n^{(1)}$. The highest weight condition is with respect to the finite
subalgebra $A_n$. The bijection $\Phi$ can be generalized by dropping the highest weight 
requirement on the elements in the KR crystals~\cite{DS:2006}, yielding the set of crystal
paths $\Path$. On the corresponding set of unrestricted rigged configurations $\RC$, the 
$A_n$ crystal structure is known explicitly~\cite{S:2006}. One of the remaining open questions is
to define the full affine crystal structure on the level of rigged configurations.
Given the affine crystal structure on both sides, the bijection $\Phi$ has a much more
conceptual interpretation as an affine crystal isomorphism.

In type $A_n^{(1)}$, the affine crystal structure can be defined using the promotion operator $\Ppr$,
which corresponds to the Dynkin diagram automorphism mapping node $i$ to $i+1$ modulo $n+1$.
On crystals, the promotion operator is defined using jeu-de-taquin~\cite{Sch:1972,Sh:2002}.
In~\cite{S:2006}, one of the authors proposed an algorithm $\Rpr$ on $\RC$ and 
conjectured~\cite[Conjecture 4.12]{S:2006} that $\Rpr$ corresponds to the promotion operator 
$\Ppr$ under the bijection $\Phi$. Several necessary conditions of promotion operators 
were established and it was shown that in special cases $\Rpr$ is the correct promotion operator.

In this paper, we show in general that $\Phi \circ \Ppr \circ \Phi^{-1}=\Rpr$ (i.e., $\Phi$ is the 
intertwiner between $\Ppr$ and $\Rpr$): 
\begin{equation*}
\begin{CD}
\Path@>{\Phi}>> \RC\\
@V{\Ppr}VV @VV{\Rpr}V \\
\Path @>>{\Phi}> \RC.
\end{CD}
\end{equation*}
Thus $\Rpr$ is indeed the promotion on $\RC$ and $\Phi$ is an affine crystal isomorphism. 

Another reformulation of the bijection from tensor products of crystals to rigged configurations 
in terms of the energy function of affine crystals and the inverse scattering formalism for the 
periodic box ball systems was given in~\cite{KS:2009,KOSTY:2006,Sa:2008,Sa:2009,Sa:2009a}.

This paper is organized as follows.
In Section~\ref{sec:definitions}, we review the definitions of crystal paths and rigged configurations,
and state the main results of this paper. Theorem~\ref{thm:main} shows that $\Rpr$ is the analogue 
of the promotion operator on rigged configurations and Corollary~\ref{cor:affine} states that $\Phi$
is an affine crystal isomorphism.
In Section~\ref{sec:outline}, we explain the outline of the proof and provide a running example
demonstrating the main ideas.
Sections~\ref{pf: Rslide and Rlh} to~\ref{pf: base case 2} contain the proofs of the results stated 
in the outline. Further technical results are delegated to the appendix.

\subsection*{Acknowledgements}
We would like to thank Nicolas Thi\'ery for his support with MuPAD-Combinat~\cite{HT:2003}
and Sage-Combinat~\cite{Sage}.
An extended abstract of this paper appeared in the FPSAC 2009 proceedings~\cite{WQ:FPSAC2009}.

%%%%%%%%%%%%%%%%%%%%%%%%%%%%%%%%%%%%%%%%%%%%%%%%%%%%%%%%%%%%%%%%%%%%%%%%%%%%%%%%%%%%%%%%%%%%%%%%%%%%%%%%%%%%%%%%%%%%%%%%%
%Section 2
%%%%%%%%%%%%%%%%%%%%%%%%%%%%%%%%%%%%%%%%%%%%%%%%%%%%%%%%%%%%%%%%%%%%%%%%%%%%%%%%%%%%%%%%%%%%%%%%%%%%%%%%%%%%%%%%%%%%%%%%%
\section{Preliminaries and the main result} \label{sec:definitions}

In this section we set up the definitions and state the main results of this paper in 
Theorem~\ref{thm:main} and Corollary~\ref{cor:affine}.
Most definitions follow~\cite{DS:2006, KSS:2002, S:2006}.

Throughout this paper the positive integer $n$ stands for the rank
of the Lie algebra $A_n$. Let $\J=[n]$ be the index set of the Dynkin diagram of type $A_n$. 
Let $\HH=\J \times \Z_{>0}$ and define $B$ to be a finite sequence of pairs of positive integers
\begin{equation*}
B = ((r_1,s_1),\ldots,(r_K,s_K)) 
\end{equation*}
with $(r_i,s_i)\in \HH$ and $1\le i \le K$.

$B$ represents a \textbf{sequence of rectangles} where the $i$-th rectangle is of height $r_i$ 
and width $s_i$. We sometimes use the phrase ``leftmost rectangle'' 
(resp.``rightmost rectangle'') to mean the first (resp. last) pair in the list. 
We use $B_i=(r_i,s_i)$ as the $i$-th pair in $B$. 

Given a sequence of rectangles $B$, we will use the following operations for successively 
removing boxes from it. In the following subsections, we define the set of paths $\Path(B)$ 
and rigged configurations $\RC(B)$, and discuss the analogous operations defined on $\Path(B)$
and $\RC(B)$. They are used to define the bijection $\Phi$ between $\Path(B)$ and $\RC(B)$ 
recursively. The proof of Theorem~\ref{thm:main} exploits this  recursion.

\begin{definition} ~\cite[Section 4.1,4.2]{DS:2006}.
\begin{enumerate}
\item If $B=((1,1), B')$, let $\lh(B)=B'$. This operation is called \textbf{left-hat}.
\item If $B=((r,s), B')$ with $s\ge 1$, let $\ls(B)=((r,1), (r,s-1), B')$.
This operation is called \textbf{left-split}. Note that when $s=1$, $\ls$ is
just the identity map.
\item If $B=((r,1), B')$ with $r\ge 2$, let $\lb(B)=((1,1),(r-1,1),B')$.
This operation is called \textbf{box-split}.
\end{enumerate}
\end{definition}

\subsection{Inhomogeneous lattice paths}
Next we define inhomogeneous lattice paths and present the analogues of the left-hat,
left-split, box-split operations on paths.

\begin{definition}
Given $(r,s)\in \HH$, define $\Path_n(r,s)$ to be 
the set of semi-standard Young tableaux of (rectangular) shape $(s^r)$ over the alphabet
$\{1,2,\ldots,n+1\}$. 
\end{definition}

Recall that for each semi-standard Young tableau $t$, we can associate a weight 
$\wt(t)=(\lambda_1, \lambda_2, \ldots, \lambda_{n+1})$ in the ambient weight lattice,
where $\lambda_i$ is the number of 
times that $i$ appears in $t$. Moreover, $\Path_n(r,s)$ is endowed with a type $A_n$-crystal 
structure, with the Kashiwara operator $e_a, f_a$ for $1\le a \le n$ defined by the signature rule. 
For a detailed discussion see for example~\cite[Chapters 7 and 8]{HK:2002}.

\begin{definition}
Given a sequence $B$ as defined above,
\begin{equation*}
	\Path_n(B) = \Path_n(r_1,s_1) \otimes \cdots \otimes \Path_n(r_K,s_K).
\end{equation*} 
\end{definition}

As a set $\Path_n(B)$ is a sequence of rectangular semi-standard Young tableaux.
It is also endowed with a crystal structure through the tensor product rule.
The Kashiwara operators $e_a,f_a$ for $1\le a\le n$ naturally extend from semi-standard 
tableaux to a list of tableaux using the signature rule. Note that in this paper we use the opposite
of Kashiwara's tensor product convention, that is, all tensor products are reverted.
For $b_1\otimes b_2 \otimes \cdots \otimes b_K \in \Path_n(B)$, 
$\wt(b_1\otimes b_2 \otimes \cdots \otimes b_K)=\wt(b_1)+\wt(b_2)+\cdots+\wt(b_K)$. 
For further details see for example~\cite[Section 2]{DS:2006}.

\begin{definition}
Let $\la=(\la_1,\la_2,\ldots,\la_{n+1})$ be a list of non-negative integers. Define
\begin{equation*}  
\Path_n(B,  \la) = \{p\in \Path_n(B) \mid \wt(p)=\la \}.
\end{equation*}  
\end{definition}

\begin{example} \label{example:path}
Let $B = ((2,2),(1,2),(3,1))$. Then
\begin{equation*}
p\;=\; \young(12,23) \otimes \young(12) \otimes \young(1,2,4) 
\end{equation*}
is an element of $\Path_3(B)$ and $\wt(p)=(3,4,1,1)$.
\end{example}

We often omit the subscript $n$, writing $\Path$ instead of $\Path_n$,
when $n$ is irrelevant or clear from the discussion. 

\begin{definition}
Let $\la=(\la_1,\la_2,\ldots,\la_{n+1})$ be a partition. Define the set of highest weight paths as
\begin{equation*}  
\Pathb_n(B, \la) = \{p\in \Path_n(B,\la) \mid \text{$e_i(p) = \emptyset$ for $i=1,2,\ldots,n$}\}.
\end{equation*}  
\end{definition}

We often refer to a rectangular tableau just as a ``rectangle'' when there is no ambiguity. For example,
the leftmost rectangle in $p$ of the above example is the tableau
\begin{equation*}
\young(12,23) \;.
\end{equation*}
  
For any $p \in \Path(B)$, the \textbf{row word} (respectively \textbf{column word}) of $p$, 
$\row(p)$ (respectively $\col(p)$), is the concatenation of the row (column) words of each rectangle 
in $p$ from left to right. 
\begin{example}
The row word of the $p$ of Example~\ref{example:path} is 
$\row(p)=\row(\young(12,23)) \cdot \row(\young(12)) \cdot \row(\young(1,2,4))=2312\cdot12\cdot421
= 231212421$, and similarly the column word is $\col(p)=213212421$.
\end{example}

\begin{definition}
We say $p \in \Path(B)$ and $q \in \Path(B^\prime)$ are Knuth equivalent, denoted by 
$p \equiv_K q$, if their row words (and hence their column words) are Knuth equivalent.
\end{definition}

\begin{example}
Let $B^\prime = ((2,2),(3,1),(1,2))$, and 
\begin{equation*}
q\;=\; \young(12,23) \otimes \young(1,2,4) \otimes \young(12) \in \Path(B^\prime)
\end{equation*}
then $p \equiv_K q$.
\end{example}

The following maps on $\Path(B)$ are the counterparts of the maps $\lh$, $\lb$ and $\ls$
defined on $B$. By abuse of notation, we use the same symbols as on rectangles.
\begin{definition}~\cite[Sections 4.1,4.2]{DS:2006}.\label{def: Pl}
\begin{enumerate}
\item Let $b=c\otimes b'\in \Path((1,1), B')$. Then $\Plh(b)=b' \in \Path(B')$.
\item Let $b=c\otimes b'\in \Path((r,s), B')$, where $c=c_1c_2\cdots c_s$ and $c_i$
denotes the $i$-th column of $c$. Then $\ls(b)=c_1\otimes c_2\cdots c_s\otimes b'$.
\item Let $b=\begin{array}{|c|} \hline b_1\\ \hline b_2\\ \hline \vdots\\ \hline b_r\\ \hline
\end{array}\otimes b'\in \Path((r,1), B')$, where $b_1<\cdots<b_r$.
Then
\begin{equation*}
	\Plb(b)=\begin{array}{|c|} \hline b_r\\ \hline \end{array} \otimes 
	\begin{array}{|c|} \hline b_1\\ \hline \vdots \\ \hline b_{r-1}\\ \hline \end{array} \otimes b'.
\end{equation*}
\end{enumerate}
\end{definition}

\subsection{Rigged configurations} \label{section.rigged_configurations}
A general definition of rigged configuration of arbitrary types can be found in
\cite[Section 3.1]{S:2006}. Here we are only concerned with type $A_n$ rigged configurations
and review their definition.

Given a sequence of rectangles $B$, following the convention of \cite{S:2006} we denote the
\textbf{multiplicity} of a given $(a,i)\in \HH$ in $B$ by setting 
$L_i^{(a)}=\#\{(r,s)\in B \mid r=a, s=i\}$.  

The (highest-weight) rigged configurations are indexed by a sequence of rectangles $B$ and a 
dominant weight $\La$. The sequence of partitions $\nu=\{\nu^{(a)}\mid a\in \J \}$ is a
\textbf{$(B,\La)$-configuration} if
\begin{equation}\label{eq:conf}
\sum_{(a,i)\in\HH} i m_i^{(a)} \alpha_a = \sum_{(a,i)\in\HH} i
L_i^{(a)} \La_a- \La,
\end{equation}
where $m_i^{(a)}$ is the number of parts of length $i$ in partition
$\nu^{(a)}$, $\alpha_a$ is the $a$-th simple root and $\La_a$ is the $a$-th 
fundamental weight. Denote the set of all $(B,\La)$-configurations by $\Conf(B,\La)$.
The \textbf{vacancy number} of a configuration is defined as
\begin{equation*}\label{eq:vac}
p_i^{(a)}=\sum_{j\ge 1} \min(i,j) L_j^{(a)}
 - \sum_{(b,j)\in \HH} (\alpha_a | \alpha_b) \min(i,j)m_j^{(b)}.
\end{equation*}
Here $(\cdot | \cdot )$ is the normalized invariant form on the weight lattice $P$
such that $A_{ab}=(\alpha_a | \alpha_b)$ is the Cartan matrix (of type $A_n$ in our case).
The $(B,\La)$-configuration $\nu$ is \textbf{admissible} if $p^{(a)}_i\ge 0$ for 
all $(a,i)\in\HH$, and the set of admissible $(B,\La)$-configurations is denoted
by $\Confb(B,\La)$.

A partition $p$ can be viewed as a linear ordering $(p,\succ)$ of a finite multiset 
of positive integers, referred to as \textbf{part}s, where parts of different lengths are ordered
by their value, and parts of the same length are given an arbitrary ordering. 
Implicitly, when we draw a Young diagram of $p$, we are giving such an ordering.
Once $\succ$ is specified, $\prec$, $\preceq$, and $\succeq$ are defined accordingly.

A labelling of a partition $p$ is then a map $J: (p,\succ) \to \Z_{\ge 0}$ satisfying
that if $i,j \in p$ are of the same value and $i \succ j$, then $J(i) \ge J(j)$ as
integers. A pair $(x, J(x))$ is referred to as a \textbf{string}, 
the part $x$ is referred to as the size or length of the string 
and $J(x)$ as its \textbf{label}.

\begin{remark} \label{rmk: succ and >}
The linear ordering $\succ$ on parts of a partition $p$ can be naturally viewed as 
an linear ordering on the corresponding strings. It is directly from its definition 
that $\succ$ is a finer ordering than $>$ that compares the size (non-negative integer) 
of the strings. Another important distinction is that $>$ can be used to compare strings
from possibly different partitions.

Given two strings $s$ and $t$, the meaning of equality $=$ is clear from the context 
in most cases. For example, if $s$ and $t$ are strings from different partitions, then $s=t$ 
means that they are of the same size; $s=t-1$ means that the length of $s$ is 
1 shorter than that of $t$. In the case that $s$ and $t$ are from the same partition and
ambiguity may arise, we reserve $s=t$ to mean $s$ and $t$ are the same string and explicitly
write $|s|=|t|$ to mean that $s$ and $t$ are of the same length but possibly distinct strings.  
\end{remark}

A \textbf{rigging} $J$ of an (admissible) $(B,\La)$-configuration $\nu = (\nu^{(1)}, \ldots, \nu^{(n)})$ 
is a sequence of maps $J=(J^{(a)})$, each $J^{(a)}$ is a labelling of the partition $\nu^{(a)}$ with the
extra requirement that for any part $i \in \nu^{(a)}$
\begin{equation*}
  0 \le J^{(a)}(i) \le p^{(a)}_i.
\end{equation*}
For each string $(i,J^{(a)}(i))$, the difference $cJ^{(a)}(i)=p_i^{(a)}-J^{(a)}(i)$ is 
referred to as the \textbf{colabel} of the string. 
$cJ=(cJ^{(a)})$ as a sequence of maps defined above is referred to as the \textbf{corigging}
of $\nu$. A string is said to be \textbf{singular} if its colabel is 0. 

\begin{definition}
The pair $\rc=(\nu,J)$ described above is called a (restricted-)\textbf{rigged configuration}.
The set of all rigged $(B,\La)$-configurations is denoted by $\RCb_n(B,\La)$.
In addition, define $\RCb(B)=\bigcup_{\La\in P^+} \RCb(B,\La)$, where $P^+$ is the set of
dominant weights.
\end{definition}

\begin{remark} \label{rmk: colabel}
Since $J$ and $cJ$ uniquely determine each other, a rigged configuration $\rc$ can be
represented either by $(\nu, J)$ or by $(\nu, cJ)$. In particular, if $x$ is a part of 
$\nu^{(a)}$ then $(x, J^{(a)}(x))$ and $(x, cJ^{(a)}(x))$ refer to the same string.
We will use these two representations interchangeably depending on which one is more 
convenient for the ongoing discussion. Nevertheless, in the later part of this paper,
when we say that a string is unchanged/preserved under some construction, we mean the 
length and the label of the string being preserved, the colabel may change due to the 
change of the vacancy number resulted from the construction.   
\end{remark}

Equation \eqref{eq:conf} provides an obvious way of defining a weight function 
on $\RCb(B)$. Namely, for $\rc \in \RCb(B)$
\begin{equation}\label{eq:weight}
\wt(\rc)=\sum_{(a,i)\in\HH} i L_i^{(a)} \La_a- \sum_{(a,i)\in\HH} i m_i^{(a)} \alpha_a .
\end{equation}

\begin{remark}\label{rmk: weight}  
When working with rigged configurations, it is often convenient to take the fundamental 
weights as basis for the weight space. On the other hand, when working with lattice paths 
we often use the ambient weight space $\Z^{n+1}$. Conceptually, this distinction is not necessary, as weights 
can be considered as abstract vectors in the weight space. One can convert from one 
representation to the other by identifying the fundamental weight $\Lambda_i$ with 
$(1^i,0^{n+1-i})$ as ambient weight. However, there is a subtlety in this conversion resulted 
from the fact that the weights are not uniquely represented by ambient weights. For example,
$(0^{n+1})$ and $(1^{n+1})$ represent the same vector in $A_n$ weight space. See 
Remark~\ref{rmk: weight again} for the conversion we use in this paper.  
\end{remark}

\begin{remark} \label{rmk: RCb(B)}
From the above definition, it is clear that $\RCb(B)$ is not sensitive to the ordering of the 
rectangles in $B$. 
\end{remark}

\begin{definition} \label{def: crystal} \cite[Section 3.2]{S:2006}   
Let $B$ be a sequence of rectangles.
Define the set of \textbf{unrestricted rigged configurations} $\RC(B)$
as the closure of $\RCb(B)$ under the operators $f_a,e_a$ 
for $a\in \J$, with $f_a,e_a$ given by:
\begin{enumerate}
\item
Define $e_a(\nu,J)$ by removing a box from a string of length $k$ in
$(\nu,J)^{(a)}$ leaving all colabels fixed and increasing the new
label by one. Here $k$ is the length of the string with the smallest
negative label of smallest length. If no such string exists,
$e_a(\nu,J)$ is undefined. 
\item
Define $f_a(\nu,J)$ by adding a box to a string of length $k$ in
$(\nu,J)^{(a)}$ leaving all colabels fixed and decreasing the new
label by one. Here $k$ is the length of the string with the smallest
non positive label of largest length. If no such string exists,
add a new string of length one and label -1.
If the result is not a valid unrestricted rigged configuration
$f_a(\nu,J)$ is undefined.
\end{enumerate}
\end{definition}

The weight function \eqref{eq:weight} defined on $\RCb(B)$ extends to $\RC(B)$ without change.

As their names suggest, $f_a$ and $e_a$ are indeed the Kashiwara operators with respect to
the weight function above, and define a crystal structure on $\RC(B)$. 
This was proved in~\cite{S:2006}.

From the definition of $f_a$, it is clear that the labels of parts in an unrestricted
rigged configuration may be negative. It is natural to ask what shapes and labels can
appear in an unrestricted rigged configuration. There is an explicit characterization 
of $\RC(B)$ which answers this question \cite[Section 3]{DS:2006}. The statement is not directly 
used in our proof, so we will just give a rough outline and
leave the interested reader to the original paper for further details: 
In the definition of $\RCb(B)$, we required that the vacancy number associated to
each part is non-negative. We dropped this requirement for $\RC(B)$. Yet the vacancy numbers   
in $\RC(B)$ still serve as the upper bound of the labels, much like the role 
a vacancy number plays for a restricted rigged configuration. For restricted rigged
configurations, the lower bound for the label of a part is uniformly 0. For unrestricted 
rigged configurations, this is not the case. The characterization gives a way on
how to find lower bound for each part.

\begin{remark} \label{rmk: RC(B)}
By Remark~\ref{rmk: RCb(B)} and Definition~\ref{def: crystal}, it is clear that $\RC(B)$ 
is not sensitive to the ordering of the rectangles in $B$. 
\end{remark}

\begin{example}
Here is an example on how we normally visualize a restricted/unrestricted 
rigged configuration. Let $B = ((2,2),(1,2),(3,1))$. Then
\begin{equation*}
	\rc = \quad \yngrc(2,-1) \quad \yngrc(1,1) \quad \yngrc(1,-1) 
\end{equation*}
is an element of $\RC(B, -\La_1+3\La_2)$.

In this example, the sequence of partitions $\nu$ is ((2),(1),(1)). 
The number that follows each part is the label assigned to this part by $J$. 
The vacancy numbers associated to these parts are $p_2^{(1)}=-1$, $p_1^{(2)}=1$, 
and $p_1^{(3)}=0$. Note that
the labels are all less than or equal to the corresponding vacancy number. In the 
case that they are equal, e.g. for the parts in $\nu^{(1)}$ and $\nu^{(2)}$, those parts
are called singular as in the case of restricted rigged configuration. 
In this example $\rc \in \RC\setminus \RCb$.
\end{example}

The following maps on $\RC(B)$ are the counterparts of $\lh$, $\lb$ and $\ls$ maps
defined on $B$. 
\begin{definition}~\cite[Section 4.1,4.2]{DS:2006} \label{def: Rl}.
\begin{enumerate}
\item Let $\rc=(\nu,J)\in\RC(B)$. Then $\Rlh(\rc) \in \RC(\lh(B))$ is defined
as follows: First set $\ell^{(0)}=1$ and then repeat the
following process for $a=1,2,\ldots,n-1$ or until stopped. Find the smallest index 
$i\ge \ell^{(a-1)}$ such that $J^{(a)}(i)$ is singular. If no such $i$ exists, set 
$\rk(\nu,J)=a$ and stop. Otherwise set $\ell^{(a)}=i$ and continue with $a+1$.
Set all undefined $\ell^{(a)}$ to $\infty$. 

The new rigged configuration $(\tilde{\nu},\tilde{J})=\Rlh(\nu,J)$ is obtained by
removing a box from the selected strings and making the new strings singular
again. 
\item Let $\rc=(\nu,J)\in\RC(B)$. Then $\Rls(\rc) \in \RC(\ls(B))$ is 
the same as $(\nu,J)$. Note however that some vacancy numbers change.
\item Let $\rc=(\nu,J)\in\RC(B)$ with $B=((r,1),B')$. Then $\Rlb(\rc) \in \RC(\lb(B))$ is 
defined by adding singular strings of length $1$ to $(\nu,J)^{(a)}$ for $1\le a < r$. 
Note that the vacancy numbers remain unchanged under $\Rlb$.
\end{enumerate}
\end{definition}

\begin{remark}
Although $\RC(B)$ does not depend on the ordering of the rectangles in $B$ 
(see Remark~\ref{rmk: RC(B)}), it is clear that the above maps depend on the ordering
in $B$.
\end{remark}

In what follows, it is often easier to work with the inverses of the above maps $\Rlh$, 
$\Rls$ and $\Rlb$ maps. In the following we give explicit descriptions of these inverses. One
can easily check that they are really inverses as their name suggests. See also~\cite{KSS:2002}.

\begin{definition} \label{def: Rl_inverse}.
\begin{enumerate}
\item 
Let $\rc \in \RC(B,\la)$ for some weight $\la$, and let $r \in [n+1]$. The map $\Rlh^{-1}$ 
takes $\rc$ and $r$ as input, and returns $\rc' \in \RC(\lh^{-1}(B),\la+\epsilon_r)$ by the 
following algorithm:   
Let $d^{(j)}=\infty$ for $j \ge r$. For $k=r-1,\ldots,1$ select the $\succ$-maximal singular string 
in $\rc^{(k)}$ of length $d^{(k)}$ (possibly of zero length) such that $d^{(k)} \le d^{(k+1)}$. 
Then $\rc'$ is obtained from $\rc$ by adding a box to each of 
the selected strings, making them singular again, and leaving
all other strings unchanged.

We denote the sequence of strings in $\rc$ selected in the above algorithm by 
\begin{equation*}
	D_r=(D^{(n)}, \ldots, D^{(1)}). 
\end{equation*}
It is called the~\textbf{$\Rlh^{-1}$-sequence} of $\rc$ 
with respect to $r$.  For simplicity for future discussions, we append $D^{(0)}=(0,0)$ to
the end of the sequence.

In light of Remark~\ref{rmk: succ and >}, we write $D^{(k)} \le D^{(k+1)}$ and say that 
$D_r$ is a weakly decreasing sequence. 

\item Let $\rc=(\nu,J)\in\RC(B)$ where $B=((r,1),(r,s),B')$. Then 
$\Rls^{-1}(\rc) \in \RC(\ls^{-1}(B))$ is the same as $(\nu,J)$.

Note that due to the change of the sequence of rectangles, the vacancy numbers
for parts in $\nu^{(r)}$ of size less than $s+1$ all decrease by 1, so the
colabels of these parts decrease accordingly. Thus $\Rls^{-1}$ is only defined
on $\rc \in \RC((r,1),(r,s),B')$ such that the colabels of parts in 
$\rc^{(k)}$ of size less than $s+1$ is $\ge 1$. All $\rc$s that satisfy the above 
conditions form $\Dom(\Rls^{-1})$. 

\item Let $\rc \in\RC(B)$ where $B=((1,1),(r-1,1),B')$. Then 
$\Rlb^{-1}(\rc) \in \RC(\lb^{-1}(B))$ is defined by removing singular strings of 
length $1$ from $\rc^{(a)}$ for $1\le a < r$, the labels of all unchanged parts 
are preserved. 

Note that the vacancy numbers remain unchanged under $\Rlb^{-1}$. As a result the 
colabels of all unchanged parts are preserved.

The collection of all $\rc \in \RC((1,1),(r-1,1),B')$ such that there is a 
singular part of size 1 in $\rc^{(a)}$ for $1 \le a < r$ forms $\Dom(\Rlb^{-1})$.
\end{enumerate}
\end{definition}

\subsection{The bijection between $\Path(B)$ and $RC(B)$} 
The map $\Phi:\Path(B,\la)\to\RC(B,\la)$ is defined recursively by various commutative 
diagrams. Note that it is possible to go from $B=((r_1,s_1),(r_2,s_2),\ldots, (r_K,s_K))$ to the empty crystal via successive application of $\lh$, $\ls$ and $\lb$. For further details 
see~\cite[Section 4]{DS:2006}.

\begin{definition} \label{def:bij}
Define the map $\Phi:\Path(B,\la)\rightarrow \RC(B,\la)$ such that 
the empty path maps to the empty rigged configuration and such that the following
conditions hold:
\begin{enumerate}
\item \label{bij:1} Suppose $B=((1,1), B')$. Then the following diagram commutes:
\begin{equation*}
\begin{CD}
\Path(B,\la) @>{\Phi}>> \RC(B,\la) \\
@V{\lh}VV @VV{\Rlh}V \\
\displaystyle{\bigcup_{\mu\in\lm} \Path(\lh(B),\mu)} @>>{\Phi}> 
\displaystyle{\bigcup_{\mu\in\lm}
\RC(\lh(B),\mu)}
\end{CD}
\end{equation*}
where $\lambda^-$ is the set of all non-negative tuples obtained from $\lambda$ by decreasing one part.
\item \label{bij:2} Suppose $B=((r,s), B')$ with $s\ge 2$. Then the 
following diagram commutes:
\begin{equation*}
\begin{CD}
\Path(B,\la) @>{\Phi}>> \RC(B,\la) \\
@V{\ls}VV @VV{\Rls}V \\
\Path(\ls(B),\la) @>>{\Phi}> \RC(\ls(B),\la).
\end{CD}
\end{equation*}
\item \label{bij:3} Suppose $B=((r,1), B')$ with $r\ge2$. Then the 
following diagram commutes:
\begin{equation*}
\begin{CD}
\Path(B,\la) @>{\Phi}>> \RC(B,\la) \\
@V{\lb}VV @VV{\Rlb}V \\
\Path(\lb(B),\la) @>>{\Phi}> \RC(\lb(B),\la).
\end{CD}
\end{equation*}
\end{enumerate}
\end{definition}

\begin{remark} \label{rmk: weight again}
By definition, $\Phi$ preserves weight. As pointed out in Remark~\ref{rmk: weight}, the
ambient weight representation is not unique. Yet for $p \in \Path(B)$, $\wt(p)$ is the 
content of $p$, which provides a ``canonical'' ambient weight representation. Passing through
$\Phi$, on $\RC(B)$ side this provides a ``canonical'' conversion between fundamental weight
and ambient weight. In particular, when we say $\rc \in \RC$ has \textbf{canonical ambient weight}
$\lambda = (\lambda_1, \ldots, \lambda_{n+1})$ we mean that $\lambda$ is the content of
$\Phi^{-1}(\rc)$. Equivalently, we are requiring that the sum of $\lambda$ is the same
as the total area of $B$
\begin{equation*} 
\sum_{i=1}^{n+1} \lambda_i = \sum_{(r,s)\in B}r\times s.
\end{equation*} 
\end{remark}

\subsection{Promotion operators}
The \textbf{promotion} operator $\Ppr$ on $\Path_n(B)$ is defined in ~\cite[page 164]{Sh:2002}.
For the purpose of our proof, we will phrase it as a composition of one \textbf{lifting} operator
and then several \textbf{sliding} operators defined on $\Path_n(B)$. 

\begin{definition}
The lifting operator $\Plift$ on $\Path_n(B)$ lifts $p \in \Path_n(B)$  to $\Plift(p) 
\in \Path_{n+1}(B)$ by adding 1 to each box in each rectangle of $p$.
\end{definition}

\begin{definition} \label{def: Pslide}
Given $p \in \Path_{n+1}(B)$, the sliding operator $\Pslide$ is defined as the 
following algorithm:
Find in $p$ the rightmost rectangle that contains $n+2$, remove one appearance of $n+2$, apply 
jeu-de-taquin on this rectangle to move the empty box to the opposite corner, fill in $1$ 
in this empty box. If no rectangle contains $n+2$, then $\Pslide$ is the identity map.

The application of jeu-de-taquin on a tableau $S$ described above naturally defines a 
\textbf{sliding route} on $S$, which is just the path along which the empty box travels from
lower right corner to upper left corner. 
\end{definition}

\begin{example}
Let $S= \young(1223,2355,4466,567)$.
After sliding lower right outside corner to the upper left inside corner, we obtain
$\Pslide(S)=\young(:123,2235,4456,5667)$. The sliding route of $S$ is \\ 
$((4,3),(3,3),(2,3),(2,2),(1,2),(1,1))$.
\end{example}

\begin{definition} \label{def: Ppr}
For $p\in \Path_n(B)$, define the promotion operator
\begin{equation*}
\Ppr(p)=\Pslide^{m}\circ\Plift(p)
\end{equation*}
where $m$ is the total number of $n+2$ in $p$. 
\end{definition}

The proposed promotion operator $\Rpr$ on $\RC_n(B)$ is defined in~\cite[Definition 4.8]{S:2006}. 
To draw the parallel with $\Ppr$ we will phrase
it as a composition of one lifting operator and then several sliding operators 
defined on $\RC(B)$. 

\begin{definition}
The lifting operator $\Rlift$ on $\RC_n(B,\la)$ lifts $\rc=(\nu,J) \in \RC_n(B,\la)$ to 
$\Rlift(\rc) \in \RC_{n+1}(B,\hat{\la})$ by setting $\Rlift(\rc)=f_1^{\lambda_1} f_2^{\lambda_2} \cdots 
f_{n+1}^{\lambda_{n+1}}(\rc)$, where 
$\lambda = (\lambda_1, \ldots, \lambda_{n+1})$ is the canonical ambient weight of $\rc$
(see Remark~\ref{rmk: weight again}) and $\hat{\la} = (0,\la_1,\ldots,\la_{n+1})$
is the canonical ambient weight of $\Rlift(\rc)$. Notice that we use the fact that
$\RC_n(B)$ is naturally embedded in $\RC_{n+1}(B)$ by simply treating the $(n+1)$-st
partition $\nu^{(n+1)}$ to be $\emptyset$.
\end{definition}

\begin{definition}\label{def: Rslide}
Given $\rc \in \RC_{n+1}(B)$, the sliding operator $\Rslide$ is defined by the following 
algorithm:
Find the $\succ$-minimal singular string in $\rc^{(n+1)}$. Let the length be $\ell^{(n+1)}$.
Repeatedly find the $\succ$-minimal singular string in $\rc^{(k)}$ of length
$\ell^{(k)}\ge \ell^{(k+1)}$ for all $1\le k<n$. Shorten the selected strings
by one and make them singular again. 

If the $\succ$-minimal singular string in $\rc^{(n+1)}$ does not exist, then 
$\Rslide$ is the identity map.

Let $I=(I^{(n+1)}, \ldots, I^{(1)}, I^{(0)})$, where for $k=n+1,\ldots,1$ the entry
$I^{(k)}$ is just the string chosen from $\rc^{(k)}$ in the above algorithm, and 
$I^{(0)}=(\infty,0)$. We call $I$ the \textbf{$\Rslide$-sequence} of $\rc$.
We say $\Rslide$ is not well-defined on $\rc$ if the $\succ$-minimal singular string
in $\rc^{(n+1)}$ exists but the $\Rslide$-sequence can not be constructed following
above algorithm (see Example~\ref{example: domain of Pslide} for what could go wrong).
\end{definition}
In light of Remark~\ref{rmk: succ and >}, we write $I^{(k)} \ge I^{(k+1)}$ and say that 
$I$ is a weakly increasing sequence. 

We note here that the above definition of $\Rslide$ is a reformulation 
of~\cite[Definition 4.8]{S:2006}.

\begin{definition} \label{def: Rpr}
Define
\begin{equation*}
\Rpr(\rc)=\Rslide^{m}\circ\Rlift(\rc)
\end{equation*}
where $m$ is the number of boxes in $\rc^{(n+1)}$.
\end{definition}

\begin{remark} \label{rmk: lifts}
It is an easy matter to show that $\Rlift=\Phi\circ\Plift\circ\Phi^{-1}$.
Indeed, we could have defined 
$\Plift(p)=f_1^{\lambda_1} f_2^{\lambda_2} \cdots f_{n+1}^{\lambda_{n+1}}(p)$, where
$\lambda = (\lambda_1, \lambda_2, \ldots, \lambda_{n+1})$ is the weight of $p$.
Since it was shown in~\cite{S:2006} that $\Phi$ is an $A_n$-crystal isomorphism, the statement
follows.
\end{remark}

There is a question in Definition~\ref{def: Ppr} on whether a sequence of $m$ $\Pslide$ 
operators can always be applied. The same question about $\Rslide$ can be asked for 
Definition~\ref{def: Rpr}. The following are examples on how things could go wrong:

\begin{example} \label{example: domain of Pslide}
Let
\begin{equation*}
p\;=\; \young(11,44) \in \Path_3(2,2). 
\end{equation*}
If we try to construct $\Pslide(p)$, we realize that after removing a copy of $4$ and move 
the empty box to the upper left corner we obtain $\young(:1,14)$, and filling the empty box 
with $1$ will violate the column-strictness of semi-standard Young tableaux. 

On the $\RC$ side, let
\begin{equation*}
\rc = \quad \emptyset \quad \yngrc(2,0) \quad \yngrc(2,0) \quad \in \RC_3(2,2).
\end{equation*}
We see that $\Rslide(\rc)$ is not well-defined. 
\end{example}

Therefore, $\Pslide$ and $\Rslide$ are partial functions on $\Path_{n+1}$ and $\RC_{n+1}$.
This, however, will not cause problems in our discussion because of the following two remarks.

\begin{remark}
$\Pslide$ is well-defined on $\Pslide^{k}(\Img(\Plift))$ for any $k$. This follows from the 
well-known fact that if $T$ is a semi-standard rectangular tableau, and if we remove all cells
that contain the largest number (which is a horizontal strip in the last row) and apply
jeu-de-taquin to move these empty cells to the upper left corner, then these empty cells
form a horizontal strip. 

$\Rslide$ is well-defined on $\Rslide^{k}(\Img(\Rlift))$ for any $k$. This is implied by~\cite[Lemma 4.10]{S:2006}.  
\end{remark}

Thus we could have just restricted the domain of $\Pslide$ to: 
\begin{definition} Define
\begin{equation*}
\Dom(\Pslide) = \bigcup_{k=0,1,2 ...} \Pslide^{k}(\Img(\Plift)).
\end{equation*}
\end{definition}

\begin{remark}
It is not known at this stage that $\Rslide$ is fully defined on $\Phi(\Dom(\Pslide))$. In fact, it 
is a consequence of our proof.
\end{remark}

Given a promotion operator in type $A_n$, we can define the affine crystal operators 
$e_0$ and $f_0$ as
\begin{equation*}
	e_0 = \Ppr^{-1} \circ e_1 \circ \Ppr \qquad \text{and} \qquad 
	f_0 = \Ppr^{-1} \circ f_1 \circ \Ppr.
\end{equation*}
An $A_n$-crystal together with $e_0$ and $f_0$ is called an \textbf{affine crystal}. An 
\textbf{affine crystal isomorphism} between crystals $B$ and $B'$ is a bijective map
$g: B \to B'$ such that $f_i\circ g(b) = g\circ f_i(b)$ for all $b\in B$ and 
$i\in\{0,1,\ldots,n\}$. See~\cite[page 164]{Sh:2002} for further discussions.

\subsection{Combinatorial $R$-matrix and right-split}
Let $B = ((r_1,s_1),\ldots,(r_K,s_K))$ be a sequence of rectangles, and let $\sigma \in S_K$ 
be a permutation of $K$ letters. $\sigma$ acts on $B$ by $\sigma(B)=((r_{\sigma(1)},s_{\sigma(1)}),
\ldots,(r_{\sigma(K)},s_{\sigma(K)}))$. 

The $R$-matrix is the affine crystal isomorphism 
$\PRmap_{\sigma}:\Path(B)\rightarrow\Path(\sigma(B))$, which sends $u_1\otimes \cdots\otimes u_K$ 
to $u_{\sigma(1)} \otimes \cdots \otimes u_{\sigma(K)}$, where $u_i\in \Path(r_i,s_i)$ is the
unique tableau of content $(s_i^{r_i})$. It was shown in~\cite[Lemma 8.5]{KSS:2002} that
for any $\sigma$, $\Phi\circ\PRmap_{\sigma}\circ\Phi^{-1} = \id$ on $\RCb(B)$. 
(Note that by Remark~\ref{rmk: RCb(B)}, $\RCb(B)$ and $\RCb(\sigma(B))$ defines the same set,
thus the above statement makes sense.)
Together with the fact that $\PRmap_{\sigma}$ preserves the $A_n$-crystal structure and
the fact that $\RC(B)$ and $\RC(\sigma(B))$ defines the same set (see Remark~\ref{rmk: RC(B)})
we have the following result.
\begin{theorem} \label{thm: PRmap and RRmap}
For any $\sigma$, $\Phi\circ\PRmap_{\sigma}\circ\Phi^{-1} = \id$ on $\RC(B)$. 
\end{theorem}

In the remainder of the paper, we often just write $\PRmap$ and omit the subscript $\sigma$. 

\begin{definition} \label{definition: rs}
$\Prs$, $\Rrs$ are called \textbf{right-split}. 
$\rs$ operates on sequences of rectangles as follows: Let $B = ((r_1,s_1),\ldots,(r_K,s_K))$, and 
suppose $s_K$ $>$ 1 (i.e, the rightmost rectangle is not a single column). Then 
$\rs(B) = ((r_1,s_1),\ldots,(r_K,s_K-1),(r_K,1))$, that is, $\rs$ splits one column off the 
rightmost rectangle.

$\Rrs$ operates on $\RC(B)$ as follows: If $\rc \in \RC(B)$, then $\Rrs(\rc) \in \RC(\rs(B))$ 
is obtained by increasing the labels by 1 for all parts in $\rc^{(r_K)}$ of size less than $s_K$. 
Observe that this will leave the colabels of all parts unchanged. 

$\Prs$, which operates on $\Path(B)$, is defined as $\Prs=\Phi\circ\Rrs\circ\Phi^{-1}$.  
\end{definition}

\subsection{The main result}
We now state the main result of this paper.
\begin{theorem}
\label{thm:main}
Let $B = ((r_1,s_1),\ldots,(r_K,s_K))$ be a sequence of rectangles, and $\Path(B)$, 
$\RC(B)$, $\Phi$, $\Ppr$, and $\Rpr$ as given as above. Then the following diagram commutes: 
\begin{equation} \label{eq:main}
\begin{CD}
\Path(B)@>{\Phi}>> \RC(B)\\
@V{\Ppr}VV @VV{\Rpr}V \\
\Path(B) @>>{\Phi}> \RC(B).
\end{CD}
\end{equation}
\end{theorem}

Using that the promotion operator on $A_n$-crystals defines an affine crystal, this also
yields the following important corollary.
\begin{corollary} \label{cor:affine}
The bijection $\Phi$ between crystal paths and rigged configurations is an affine crystal
isomorphism.
\end{corollary}

\section{Outline of the proof of Theorem~\ref{thm:main}} \label{sec:outline}
In this section, we draw the outline of the proof and state all important results needed in the proof, 
but leave the details of the proofs to later sections. We also illustrate the main ideas
with a running example.

By Remark~\ref{rmk: lifts}, for the proof of Theorem~\ref{thm:main} it suffices to show 
that the following diagram commutes:
\begin{equation*}
\begin{CD}
\Dom(\Pslide)@>{\Phi}>> \Phi(\Dom(\Pslide))\\
@V{\Pslide}VV @VV{\Rslide}V \\
\Dom(\Pslide) @>>{\Phi}> \Phi(\Dom(\Pslide)).
\end{CD}
\end{equation*}
In particular, we need to show that $\Rslide$ is defined on $\Phi(\Dom(\Pslide))$.

\subsection{Setup the running example}
As an abbreviation, for any $p \in \Dom(\Pslide)$, we use $\D(p)$ to mean the following statement:
 ``$\Rslide(\Phi(p))$ is well-defined and the diagram
\begin{equation*}
\begin{CD}
p@>{\Phi}>> \bullet \\ 
@V{\Pslide}VV @VV{\Rslide}V \\
\bullet @>>{\Phi}> \bullet
\end{CD}
\end{equation*}
commutes''.

For $p,q \in \Dom(\Pslide)$ we write $\D(p)\leadsto\D(q)$ to mean that $\D(p)$ reduces to $\D(q)$,
that is, $\D(q)$ is a sufficient condition for $\D(p)$.

We will let $n=3$ and use the following $p \in \Path_3((2,2),(3,2),(2,2))$ as the starting point 
of the running example:   
\begin{equation*}
p\;=\; \young(22,44) \otimes \young(12,23,34) \otimes \young(12,23).
\end{equation*}

After lifting to $\Path_4$ we have:
\begin{equation*}
l(p)\;=\; \young(33,55) \otimes \young(23,34,45) \otimes \young(23,34) \in \Dom(\Pslide).
\end{equation*}

Our goal is to show $\D(l(p))$ by a sequence of reductions.
Note that the rightmost $5$ (which is $n+2$ for $n$=3) appears in the second rectangle. Thus
$\Pslide$ acts on the second rectangle. The first motivation behind our reductions is to try 
to get rid of boxes from the left and make $\Pslide$ act on the leftmost rectangle:
\begin{description} 
\item [Step 1]
\begin{equation*}
\D(\young(33,55) \otimes \young(23,34,45) \otimes \young(23,34)) 
\stackrel{\ls}{\leadsto}
\D(\young(3,5) \otimes \young(3,5) \otimes \young(23,34,45) \otimes \young(23,34))
\end{equation*}
This is called a $\ls$-reduction, which is justified by Theorems~\ref{thm: Pslide and Pls} 
and~\ref{thm: Rslide and Rls} below.

\item [Step 2]
\begin{equation*}
\D(\young(3,5) \otimes \young(3,5) \otimes \young(23,34,45) \otimes \young(23,34))
\stackrel{\lb}{\leadsto}
\D(\young(5) \otimes \young(3) \otimes \young(3,5) \otimes \young(23,34,45) \otimes \young(23,34))
\end{equation*}
This is called a $\lb$-reduction, which is justified by Theorems~\ref{thm: Pslide and Plb} 
and~\ref{thm: Rslide and Rlb} below.

\item [Step 3]
\begin{equation*}
\D(\young(5) \otimes \young(3) \otimes \young(3,5) \otimes \young(23,34,45) \otimes \young(23,34))
\stackrel{\lh}{\leadsto}
\D(\young(3) \otimes \young(3,5) \otimes \young(23,34,45) \otimes \young(23,34))
\end{equation*}
This is called a $\lh$-reduction, which is justified by Theorems~\ref{thm: Pslide and Plh} 
and~\ref{thm: Rslide and Rlh} below.

\item [Step 4]
Another application of $\lh$-reduction.
\begin{equation*}
\D(\young(3) \otimes \young(3,5) \otimes \young(23,34,45) \otimes \young(23,34))
\stackrel{\lh}{\leadsto}
\D(\young(3,5) \otimes \young(23,34,45) \otimes \young(23,34))
\end{equation*}

\end{description} 

We repeat above reductions until the rightmost tableau containing 5 becomes the first tableau
in the list. After that we want to further simplify the list, if possible, to get rid of boxes 
from right by pushing them column-by-column to the left using the $\PRmap$-matrix map $\PRmap$, 
until we reach the place where can prove $\D(\bullet)$ directly: 

\begin{description} 
\item [Step 8]

\begin{equation*}
\D(\young(23,34,45) \otimes \young(23,34))
\stackrel{\rs}{\leadsto}
\D(\young(23,34,45) \otimes \young(3,4) \otimes \young(2,3))
\end{equation*}
This is called a $\rs$-reduction, which is justified by Theorems~\ref{thm: Pslide and Prs} 
and~\ref{thm: Rslide and Rrs} below.

\item [Step 9]
\begin{equation*}
\D(\young(23,34,45) \otimes \young(3,4) \otimes \young(2,3))
\stackrel{\PRmap}{\leadsto}
\D(\young(3,4) \otimes \young(23,34,45) \otimes \young(2,3))
\end{equation*}
This is called a $\PRmap$-reduction, which is justified by Theorem~\ref{thm: Pslide and PRmap}.
\end{description} 

Now since the rectangle that $\Pslide$ acts on is no longer the leftmost one, we can go back to Step 1.  
Repeating the above steps until $\Pslide$ acts on the leftmost rectangle again, we need one more
$\PRmap$-reduction: 

\begin{description} 
\item [Step 13]
\begin{equation*}
\D(\young(23,34,45) \otimes \young(2,3))
\stackrel{\PRmap}{\leadsto}
\D(\young(3,5) \otimes \young(22,33,44))
\end{equation*}
\end{description} 

Using these reductions, we will eventually reach one of the following two base cases:
\begin{itemize}
  \item Base case 1: $p$ is a single rectangle that contains $n+2$; or
  \item Base case 2: $p= S \otimes q$, where $S$ is a single column that contains $n+2$, 
  and $n+2$ does not appear in $q$. 
\end{itemize}
In certain cases it might be possible to reduce Base case 2 further to Base case 1. But we 
will prove both base cases in this full generality without specifying when this further reduction 
is possible.

In the above example, we reached the second case. 
\begin{description} 
\item [Step 14]
Now we have to prove this base case directly:
\begin{equation*}
\D(\young(3,5) \otimes \young(22,33,44))
\end{equation*}
This is justified by Theorem~\ref{thm: base case 2}. Base case 1 is proved in
Theorem~\ref{thm: base case 1}.
\end{description} 

\subsection{The reduction}
In this section, we formalize the ideas demonstrated in the previous section.

\begin{definition}  
Define
\begin{equation*}
\PBase=\{ p \in \Dom(\Pslide) \mid \text{$n+2$, if any exist, appears only in the leftmost rectangle of $p$}\}.
\end{equation*}
\end{definition}

The next two theorems concern the $\lh$-reduction: $\D(p) \stackrel{\lh}{\leadsto} \D(\lh(p))$.
 
\begin{theorem} \label{thm: Pslide and Plh}
   Let $p \in (\Dom(\Pslide)\setminus\PBase) \cap \Dom(\Plh)$. Then $\Plh(p) \in \Dom(\Pslide)$
   and the following diagram commutes:
   \begin{equation*}
   \begin{CD}
   p@>{\Plh} >> \bullet \\ 
   @V{\Pslide}VV @VV{\Pslide}V \\
   \bullet @>>{\Plh} > \bullet
   \end{CD}
   \end{equation*}
   \end{theorem}
   \begin{proof}
   By definition, $\Pslide$ acts on the rightmost rectangle of $p$ that contains the number $n+2$.
   Given $p \in (\Dom(\Pslide)\setminus\PBase) \cap \Dom(\Plh)$, the rightmost rectangle that contains
   $n+2$ is not the leftmost one in $p$, thus $\Pslide$ does not act on the leftmost rectangle of $p$.
   But $\Plh$, by definition, acts on the leftmost rectangle, and it is clear that 
   $\Plh(p) \in \Dom(\Pslide)$ if $p \in \Dom(\Pslide)$, and that the diagram commutes.
   \end{proof}

\begin{theorem} \label{thm: Rslide and Rlh}
   Let $\rc \in \Phi((\Dom(\Pslide)\setminus\PBase) \cap \Dom(\Plh))$ and assume that 
   $\Rslide(\Rlh(\rc))$ is
   well-defined. Then $\Rslide(\rc)$ is well-defined and the following diagram commutes: 
   \begin{equation*}
   \begin{CD}
   \rc@>{\Rlh} >> \bullet \\ 
   @V{\Rslide}VV @VV{\Rslide}V \\
   \bullet @>>{\Rlh}> \bullet
   \end{CD}
   \end{equation*}
   \end{theorem}
   \begin{proof}
   See Section~\ref{pf: Rslide and Rlh}.
   \end{proof}

To see that the above two theorems suffice for the $\lh$-reduction, we let $p$ and $\rc$ be given as 
above and consider the following diagram 
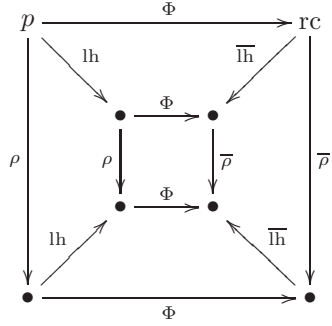
\begin{equation*}
\xymatrix{
 {p} \ar[rrr]^{\Phi} \ar[ddd]_{\Pslide} \ar[dr] ^{\Plh}& & &
        {\rc} \ar[ddd]^{\Rslide} \ar[dl]_{\Rlh} \\
 & {\bullet} \ar[r]^{\Phi} \ar[d]_{\Pslide} & {\bullet} \ar[d] ^{\Rslide}& \\
 & {\bullet} \ar[r]^{\Phi}  & {\bullet}  & \\
 {\bullet} \ar[rrr]_{\Phi} \ar[ur] ^{\Plh}& & & {\bullet} \ar[ul]_{\Rlh}
}
\end{equation*}
This diagram should be viewed as a ``cube'', the large outside square being the front face and
the small inside square being the back face, the four trapezoids between these two squares
are the upper, lower, left and right faces, respectively. We observe the following:
\begin{enumerate} 
  \item The upper and lower face commute by \cite{DS:2006}.
  \item By Theorem \ref{thm: Pslide and Plh}, the left face commutes.
  \item If we assume that the back face commutes, in particular that $\Rslide$ on the right
        edge of the back face is well-defined, then by Theorem~\ref{thm: Rslide and Rlh}
        we can conclude that the right face is well-defined and commutes.
\end{enumerate} 
Thus, if we assume the commutativity of the back face, the commutativity of the front 
face follows by induction.

The next two theorems are for $\lb$-reduction: $\D(p) \stackrel{\lb}{\leadsto} \D(\lb(p))$.
\begin{theorem} \label{thm: Pslide and Plb}
   Let $p \in (\Dom(\Pslide)\setminus\PBase) \cap \Dom(\Plb)$. Then $\Plb(p) \in \Dom(\Pslide)$ 
   and the following diagram commutes:
   \begin{equation*}
   \begin{CD}
   p@>{\lb}>> \bullet \\ 
   @V{\Pslide}VV @VV{\Pslide}V \\
   \bullet @>>{\lb}> \bullet
   \end{CD}
   \end{equation*}
\end{theorem}
\begin{proof}
The proof is similar to the argument for the $\lh$-reduction (see Theorem~\ref{thm: Pslide and Plh}).
\end{proof}

\begin{theorem} \label{thm: Rslide and Rlb}
   Let $\rc \in \Phi((\Dom(\Pslide)\setminus\PBase) \cap \Dom(\Plb))$ and assume that 
   $\Rslide(\Rlb(\rc))$ is well-defined. Then both $\Rslide(\rc)$ and $\Rlb(\Rslide(\rc))$
   are well-defined and the following diagram commutes:
   \begin{equation*}
   \begin{CD}
   \rc@>{\Rlb}>> \bullet \\ 
   @V{\Rslide}VV @VV{\Rslide}V \\
   \bullet @>>{\Rlb}> \bullet
   \end{CD}
   \end{equation*}
   \end{theorem}
   \begin{proof}
   See Section~\ref{pf: Rslide and Rlb}.
   \end{proof}

The reason that the above two theorems suffice for the $\lb$-reduction is analogous to the reason for the 
$\lh$-reduction.
 
The next two theorems are for $\ls$-reduction: $\D(p) \stackrel{\ls}{\leadsto} \D(\ls(p))$.
\begin{theorem} \label{thm: Pslide and Pls}
   Let $p \in (\Dom(\Pslide)\setminus\PBase) \cap \Dom(\Pls)$. Then $\Pls(p) \in \Dom(\Pslide)$ 
   and the following diagram commutes:
   \begin{equation*}
   \begin{CD}
   p@>{\ls}>> \bullet \\ 
   @V{\Pslide}VV @VV{\Pslide}V \\
   \bullet @>>{\ls}> \bullet
   \end{CD}
   \end{equation*}
\end{theorem}
\begin{proof}
The proof is similar to the argument for $\lh$-reduction (see Theorem~\ref{thm: Pslide and Plh}).
\end{proof}

\begin{theorem} \label{thm: Rslide and Rls}
   Let $\rc \in \Phi((\Dom(\Pslide)\setminus\PBase) \cap \Dom(\Pls))$ and assume that 
   $\Rslide(\Rls(\rc))$ is well-defined. Then $\Rslide(\rc)$ is well-defined and the following 
   diagram commutes:
   \begin{equation*}
   \begin{CD}
   \rc@>{\Rls}>> \bullet \\ 
   @V{\Rslide}VV @VV{\Rslide}V \\
   \bullet @>>{\Rls}> \bullet
   \end{CD}
   \end{equation*}
   \end{theorem}
   \begin{proof}
   See Section~\ref{pf: Rslide and Rls}.
   \end{proof}

The reason that the above two theorems suffice for the $\ls$-reduction is analogous to the reason for the 
$\lh$-reduction.

The above $\lh/\lb/\ls$-reductions 
make it clear that we have $\D(p)$
for any $p \in (\Dom(\Pslide)\setminus\PBase)$, thus reducing the problem to proving $\D(p)$ for
$p \in \PBase$. 

For $p \in \PBase$, $\D(p)$ is proved by another round of reductions, until $p$ is in 
one of the two base cases (not mutually exclusive): 

\begin{definition} [Base case 1] \label{definition: base case 1}
\begin{equation*}
\PBBaseA = \{p \in \PBase \mid \text{$p$ is a single rectangle} \} .
\end{equation*}
\end{definition}

\begin{definition} [Base case 2] \label{definition: base case 2}
\begin{equation*}
\PBBaseB = \{p \in \PBase \mid \text{the leftmost rectangle of $p$ is a single column} \} .
\end{equation*}
\end{definition}

The next two theorems deal with $\rs$-reduction: $\D(p) \stackrel{\rs}{\leadsto} \D(\rs(p))$.
\begin{theorem} \label{thm: Pslide and Prs}
   Let $p \in \PBase \setminus \PBBaseA$. Then $\Prs(p) \in \PBase \setminus \PBBaseA$
   and the following diagram commutes:
   \begin{equation*}
   \begin{CD}
   p@>{\Prs}>> \bullet \\ 
   @V{\Pslide}VV @VV{\Pslide}V \\
   \bullet @>>{\Prs}> \bullet
   \end{CD}
   \end{equation*}
\end{theorem}
\begin{proof}
See Section~\ref{pf: Pslide and Prs}.
\end{proof}

\begin{theorem} \label{thm: Rslide and Rrs}
   Let $\rc \in \Phi(\PBase \setminus \PBBaseA)$ and assume that $\Rslide(\Rrs(\rc))$ is well-defined.
   Then $\Rslide(\rc)$ is well-defined and the following diagram commutes:
   \begin{equation*}
   \begin{CD}
   \rc@>{\Rrs}>> \bullet \\ 
   @V{\Rslide}VV @VV{\Rslide}V \\
   \bullet @>>{\Rrs}> \bullet
   \end{CD}
   \end{equation*}
\end{theorem}
\begin{proof}
By the definition of $\Rslide$ and by the fact that $\Rrs$ preserves the colabels of all parts, 
it is clear that if the $\Rslide$-sequence of $\Rrs(\rc)$ exists, then the $\Rslide$-sequence of $\rc$ 
must exist and be the same as that for $\Rrs(\rc)$. Then commutativity follows.
\end{proof}

To see that the above two theorems suffice for the $\rs$-reduction, we let $p$ and $\rc$ be given as 
above and consider the following diagram 
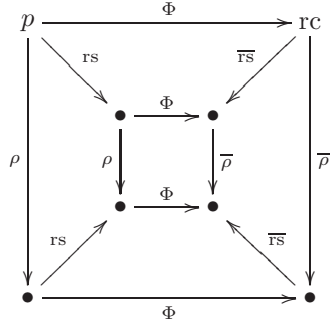
\begin{equation*}
\xymatrix{
 {p} \ar[rrr]^{\Phi} \ar[ddd]_{\Pslide} \ar[dr] ^{\Prs}& & &
        {\rc} \ar[ddd]^{\Rslide} \ar[dl]_{\Rrs} \\
 & {\bullet} \ar[r]^{\Phi} \ar[d]_{\Pslide} & {\bullet} \ar[d] ^{\Rslide}& \\
 & {\bullet} \ar[r]^{\Phi}  & {\bullet}  & \\
 {\bullet} \ar[rrr]_{\Phi} \ar[ur] ^{\Prs}& & & {\bullet} \ar[ul]_{\Rrs}
}
\end{equation*}

We observe the following:
\begin{enumerate} 
  \item The upper and lower face commute by the definition of $\Rrs$ and $\Prs$ as stated
  in Definition~\ref{definition: rs}.
  \item The left face commutes by Theorem~\ref{thm: Pslide and Prs}.
  \item If we assume that the back face commutes, in particular $\Rslide$ on the right
        edge of the back face is well-defined, then by Theorem~\ref{thm: Rslide and Rrs}
        we can conclude that the right face is well-defined and commutes.
\end{enumerate} 
Thus, if we assume the commutativity of the back face, the commutativity of the front 
face follows. 

The next theorem is for $\PRmap$-reduction: $\D(p) \stackrel{\PRmap}{\leadsto} \D(\PRmap(p))$.
\begin{theorem} \label{thm: Pslide and PRmap}
   Let $p \in \PBase \subset \Path(B)$ where $B=((r_1,s_1),(r_2,s_2))$.
   Then $\PRmap(p) \in \Dom(\Pslide)$ and the following diagram commutes:
   \begin{equation*}
   \begin{CD}
   p@>{\PRmap}>> \bullet \\ 
   @V{\Pslide}VV @VV{\Pslide}V \\
   \bullet @>>{\PRmap}> \bullet
   \end{CD}
   \end{equation*}
\end{theorem}
\begin{proof}
   It was shown in~\cite[Lemma 5.5, Eq. (5.8)]{SW:1999} that $\PRmap$ and
   $\Pslide$ commute on standardized highest weight paths (the maps are called $\sigma_i$ 
   and $\mathcal{C}_p$ in~\cite{SW:1999}, respectively). For a given $p$, we can always find a 
   $q\in \Path(B')$ for some $B'$ such that $p\otimes q$ is highest weight and $q$ does not contain 
   any $n+2$ (basically $q$ needs to be chosen such that $\varphi_i(q)\ge \varepsilon_i(p)$ for 
   all $i=1,2,\ldots,n+1$). Since $\PRmap$ respects Knuth relations, it is well-behaved with respect
   to standardization. Similarly, $\rho$ is well-behaved with respect to standardization because
   jeu-de-taquin is.
   Since by assumption $p\in\Dom(\rho)$, this implies the statement of the theorem.
\end{proof}

As the next remark shows, we only need Theorem~\ref{thm: Pslide and PRmap} in the special
case $s_2=1$. An independent proof of Theorem~\ref{thm: Pslide and PRmap} for $s_2=1$
will appear in the PhD thesis of the second author.

\begin{remark} \label{rmk: Rmap}
We would like to point out that Theorem~\ref{thm: Pslide and PRmap} for $s_2=1$
suffices for the $\PRmap$-reduction.

By definition, $\Pslide$ acts on the rightmost rectangle that contains $n+2$. If $p \in \PBase$,
then the rightmost rectangle that contains $n+2$ is also the leftmost (thus the only) rectangle
that contains $n+2$. This implies that if some permutation $\sigma$ does not involve swapping
the first two rectangle (that is, $s_1$ does not appear in the reduced word of $\sigma$), then
$\Pslide$ clearly commutes with $\PRmap_{\sigma}$.

Without loss of generality, we can further assume that the second rectangle is a single 
column. For if it is not, we can use right-split to split off a single column from the rightmost rectangle
(which commutes with $\Pslide$ by Theorem~\ref{thm: Rslide and Rrs}). Then we can use the
$\PRmap$ to move this single column to be the second rectangle (which commutes with 
$\Pslide$ by above argument). Hence it suffices to consider the case $B=((r_1,s_1),(r_2,1))$.

It is worth pointing out that although Theorem~\ref{thm: Pslide and PRmap} only
states the commutativity of $\Pslide$ and $\PRmap$ in this special case, that as a consequence 
of our main result Theorem~\ref{thm:main}, $\Pslide$ and $\PRmap$ commute in general.
\end{remark}

To see that the above theorem suffices for the $\PRmap$-reduction, we consider the 
following diagram 
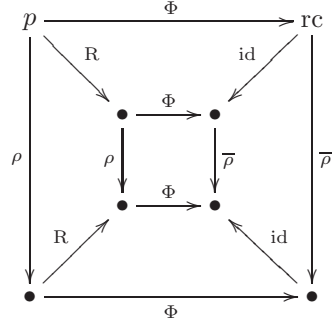
\begin{equation*}
\xymatrix{
 {p} \ar[rrr]^{\Phi} \ar[ddd]_{\Pslide} \ar[dr] ^{\PRmap}& & &
        {\rc} \ar[ddd]^{\Rslide} \ar[dl]_{\id} \\
 & {\bullet} \ar[r]^{\Phi} \ar[d]_{\Pslide} & {\bullet} \ar[d] ^{\Rslide}& \\
 & {\bullet} \ar[r]^{\Phi}  & {\bullet}  & \\
 {\bullet} \ar[rrr]_{\Phi} \ar[ur] ^{\PRmap}& & & {\bullet} \ar[ul]_{\id}
}
\end{equation*}

We observe the following:
\begin{enumerate} 
  \item The upper and lower face commute by definition of $\PRmap$. 
  \item The left face commutes by Theorem~\ref{thm: Pslide and PRmap}.
  \item The right face commutes trivially. 
\end{enumerate} 
Thus, if we assume the commutativity of the back face, the commutativity of the front 
face follows. 

Finally, we state the theorems for dealing with the base cases:
\begin{theorem}[Base case 1] \label{thm: base case 1}
Let $p \in \PBBaseA$, then $\D(p)$.
\end{theorem}
\begin{proof}
See Section~\ref{pf: base case 1}.
\end{proof}

\begin{theorem}[Base case 2] \label{thm: base case 2}
Let $p \in \PBBaseB$, then $\D(p)$. 
\end{theorem}
\begin{proof}
See Section~\ref{pf: base case 2}.
\end{proof}

%%%%%%%%%%%%%%%%%%%%%%%%%%%%%%%%%%%%%%%%%%%%%%%%%%%%%%%%%%%%%%%%%%%%%%%%%%%%%%%%%%%%%%%%%%%%%%%%%%%%%%%%%%%%%%%%%%%%%%%%%
%
%%%%%%%%%%%%%%%%%%%%%%%%%%%%%%%%%%%%%%%%%%%%%%%%%%%%%%%%%%%%%%%%%%%%%%%%%%%%%%%%%%%%%%%%%%%%%%%%%%%%%%%%%%%%%%%%%%%%%%%%%
\section{Proof of Theorem~\ref{thm: Rslide and Rlh}} \label{pf: Rslide and Rlh} 

The statement of Theorem~\ref{thm: Rslide and Rlh} is clearly equivalent to the 
following statement: 
Let $\rc \in \RC$ be such that $\Rslide$ is well-defined on $\rc$ and 
$\Rlh^{-1}(\rc,r) \in \Phi((\Dom(\Pslide)\setminus\PBase)$ for some $r \in [n+2]$.
Then $\Rslide$ is well-defined on $\Rlh^{-1}(\rc,r)$ and the following diagram commutes:
   \begin{equation*}
   \begin{CD}
   \rc@>{\Rlh^{-1}} >> \bullet \\ 
   @V{\Rslide}VV @VV{\Rslide}V \\
   \bullet @>>{\Rlh^{-1}}> \bullet
   \end{CD}
   \end{equation*}
Indeed, this is the statement we are going to prove.

Let us first consider the case that $\Rslide$ is the identity map on $\rc$. The map $\Rslide$ being the 
identity means that $\rc^{(n+1)}$ does not contain any singular string. If $r<n+2$, then 
$\Rlh^{-1}(\rc,r)^{(n+1)}$ still does not contain any singular string since no strings or vacancy
numbers in the $(n+1)$-st rigged partition change. Thus $\Rslide$ is the identity map on 
$\Rlh^{-1}(\rc,r)$. Clearly $\Rlh^{-1}$ and $\Rslide$ commute.

If $r=n+2$, then by Definition~\ref{def: Rl_inverse}, the $\Rlh^{-1}$-sequence of $\rc$ with 
respect to $n+1$ is a sequence of all 0s. Thus for each $k$, $\Rlh^{-1}(\rc)^{(k)}$ has a 
singular sting of size 1. Therefore the $\Rslide$-sequence of $\Rlh^{-1}(\rc)$ exists and
is a sequence of all 1s. Combining with the fact that $\Rlh^{-1}$ and $\Rslide$ preserve
all unchanged strings we can conclude that $\Rlh^{-1}$ and $\Rslide$ commute.

From now on, we shall assume that $\Rslide$ is not the identity map on $\rc$

Let $D_r$ be the $\Rlh^{-1}$-sequence given in Definition~\ref{def: Rl_inverse}.
Let $I$ be the $\Rslide$-sequence given in Definition~\ref{def: Rslide}.
We note that by definition $I^{(n+1)} \not\succ D_r^{(n+1)}$ and $I^{(0)} \succ D_r^{(0)}$. 
Thus, one of the following two statements must hold:
\begin{enumerate} 
  \item There is an index $N \in \{1,\ldots,n+1\}$ such that $D_r^{(N)} \succ I^{(N)}$ and 
  $D_r^{(N-1)} \prec I^{(N-1)}$;
  \item There is an index $N \in \{1,\ldots,n+1\}$ such that $D_r^{(N)} = I^{(N)}$.
\end{enumerate} 

\begin{remark} \label{rmk: cross point}
In either case above, we say $D_r$ and $I$ cross at the position $N$.
\end{remark}

Let $\rc=(u,U)$, $\Rlh^{-1}(\rc)=(v,V)$, $\Rslide(\rc)=(w,W)$, and 
$\Rlh^{-1}\circ\Rslide(\rc)=(x,X)$, $\Rslide\circ\Rlh^{-1}(\rc)=(y,Y)$. 
We denote by $\Rslide(D)$ the $\Rlh^{-1}$-sequence of $\Rslide(\rc)$, 
and denote by $\Plh^{-1}(I)$ the $\Rslide$-sequence of $\Plh^{-1}(\rc)$.

The readers may want to review Remark~\ref{rmk: succ and >} for notations used in 
the following proof.

\subsection{Case 1}
In this case we must have $I^{(N-1)} > D^{(N)}$ and $D^{(N-1)} < I^{(N)}$.
This then implies that $D^{(N-1)}<I^{(N-1)}-1$, from which we can conclude (considering
the changes in vacancy numbers) that 
\begin{itemize}
   \item $\Rslide(D)^{(k)}=D^{(k)}$ for $k \neq N$, and $\Rslide(D)^{(N)}=I^{(N)}-1$;
   \item $\Rlh^{(-1)}(I)^{(k)}=I^{(k)}$ for $k \neq N$, and $\Rlh^{-1}(I)^{(N)}=D^{(N)}+1$.
\end{itemize}

To show $(x,X)=(y,Y)$, we first argue that $x^{(k)} = y^{(k)}$, then we  
show that on corresponding parts of $x^{(k)}$ and $y^{(k)}$, 
$X^{(k)}$ and $Y^{(k)}$ either agree on their labels or agree on their colabels.
All the above and the fact that $(x,X)$ and $(y,Y)$ have the 
same sequence of rectangles implies that $X=Y$ thus $(x,X)=(y,Y)$. 

We divide the argument into the following three cases: 
\begin{itemize}
  \item $k > N$
  \item $k=N$;
  \item $1 \le k < N$;
\end{itemize}

For $k > N$, since $\Rslide(D)^{(k)}=D^{(k)}$ and $\Rlh^{-1}(I)^{(k)}=I^{(k)}$ we 
know $x^{(k)}=y^{(k)}$, both differ from $u^{(k)}$ by ``moving a box from $I^{(k)}$ to $D^{(k)}$''.
Furthermore, by the definition of $\Rslide$ and $\Rlh^{-1}$, in both $x^{(k)}$ and $y^{(k)}$, 
the labels of all unchanged strings are preserved. In the two changed strings, 
one gets a box removed and one gets a box added, and they are both kept singular.

For $k=N$, $\Rlh^{-1}(I)^{(N)}=D^{(N)}+1$ implies that from $u^{(N)}$ to $v^{(N)}$ to $y^{(N)}$ 
one box is added to $D^{(N)}$ and then is removed, thus keeping $y^{(N)}=u^{(N)}$.
Similarly, $\Rslide(D)^{(N)}=I^{(N)}-1$ implies that from $u^{(N)}$ to $w^{(N)}$ to $x^{(N)}$ 
one box is removed from $I^{(N)}$ and then is added back, thus keeping $x^{(N)}=u^{(N)}$. 
Hence $x^{(N)}=y^{(N)}$. From $U^{(N)}$ to $V^{(N)}$ to $Y^{(N)}$ the labels of all strings 
other than $D^{(N)}$ are unchanged, and for part $D^{(N)}$, both $\Rlh^{-1}$ and $\Rslide$ 
preserve its singularity. 

For $1 \le k < N$, a similar argument as the $k > N$ shows the desired result.  

\subsection {Case 2}
Let $N=\max\{k \mid I^{(k)} = D^{(k)}\}$ and let $M=\max\{k \mid I^{(k)} > D^{(k)}\}$. Thus
clearly $M < N$ and for $k>N$, $D^{(k)}\succ I^{(k)}$; for $M < k \le N$, 
$|I^{(k)}|=|D^{(k)}|$ (it may not be the case that $I^{(k)} = D^{(k)}$); 
for $k \le M$, $I^{(k)} \succ D^{(k)}$, 
in particular for $k=M$, $D^{(M)}<D^{(M+1)}$ and $I^{(M+1)}<I^{(M)}$, so 
$D^{(k)}<I^{(k)}-1$ for $k \le M$.

The above discussion implies that
\begin{itemize}
   \item $\Rslide(D)^{(k)}=D^{(k)}$ and $\Rlh^{-1}(I)^{(k)}=I^{(k)}$ for $k > N$;
   \item $\Rslide(D)^{(k)}=I^{(k)}-1$ and $\Rlh^{-1}(I)^{(k)}=D^{(k)}+1$ for $M < k \le N$;
   \item $\Rslide(D)^{(k)}=D^{(k)}$ and $\Rlh^{-1}(I)^{(k)}=I^{(k)}$ for $k \le M$.
\end{itemize}

Following the same strategy as in Case 1, we divide our argument into the
following three cases: 
\begin{itemize}
  \item $k > N$;
  \item $M < k \le N$;
  \item $1 \le k \le M$.
\end{itemize}
 
For $k > N$, the argument is the same as the $k > N$ discussion of Case 1. 

For $M <k \le N$, $\Rlh^{-1}(I)^{(k)}=D^{(k)}+1$ implies that from $u^{(k)}$ to $v^{(k)}$ to 
$y^{(k)}$ one box is added to $D^{(k)}$ and then is removed, thus keeping $y^{(k)}=u^{(k)}$.
Similarly, $\Rslide(D)^{(k)}=I^{(k)}-1$ implies that from $u^{(k)}$ to $w^{(k)}$ to $x^{(k)}$ 
one box is removed from $I^{(k)}$ and then is added back, thus keeping $x^{(k)}=u^{(k)}$. 
Hence $x^{(k)}=y^{(k)}$. From $U^{(k)}$ to $V^{(k)}$ to $Y^{(k)}$ the labels of all parts 
other than $D^{(k)}$ are unchanged, and for part $D^{(k)}$, both $\Rlh^{-1}$ and $\Rslide$ 
preserve its singularity. Moreover, the vacancy number of parts of size $|D^{(k)}|$ 
is unchanged from $U^{(k)}$ to $Y^{(k)}$ due to the cancellation of the effects of
removing $D^{(k-1)}$ (or changing the sequence of rectangles for the case $N=1$) and adding
$I^{(k+1)}$. Thus the label of $D^{(k)}$ is unchanged from $U^{(k)}$ to $Y^{(k)}$. Hence 
$U^{(k)}=Y^{(k)}$. An analogous argument shows that $U^{(k)}=X^{(k)}$. Thus $X^{(k)}=Y^{(k)}$.

For $k \le M$, the argument is similar to the case $k>N$. 

\subsection{Some remark}
We could have in both cases above defined $M=\max\{k<N \mid I^{(k)} > I^{(k+1)}\}$. Then it
would agree with the $M$ defined in Case 2, and the proof of Case 2 
could conceptually unify the two cases into one argument, but it probably would not
make the proof more readable. But this definition of $M$ does simplify statement like the 
following. 

\begin{lemma} \label{lemma: I increases under Rlh inverse}
For any $k \in [n+1]$, $\Rlh^{-1}(I)^{(k)} \ge I^{(k)}$. 
The strict inequality $\Rlh^{-1}(I)^{(k)} < I^{(k)}$ is obtained precisely on 
$(M,N]$. In particular, $\Rlh^{-1}(I)^{(k)}=I^{(k)}$ if $D^{(k)}>I^{(k)}$.
\end{lemma}
The lemma follows from the proof in Case 1 and 2, and will be referred to in the future
sections. 

\begin{remark} \label{rmk: converse}
The same idea used in the proof of this section can be used to prove the following converse of 
Theorem~\ref{thm: Rslide and Rlh}, which will be used in the proof of 
Theorem~\ref{thm: Rslide and Rlb} in Section~\ref{pf: Rslide and Rlb}.
 
\begin{prop} \label{prop: Rslide and Rlh}
Let $\rc \in \RC$ be such that $\Rslide(\rc)$ is well-defined.
Then $\Rslide$ is well-defined on $\Rlh(\rc)$ and the following diagram commutes:
   \begin{equation*}
   \begin{CD}
   \rc@>{\Rlh} >> \bullet \\ 
   @V{\Rslide}VV @VV{\Rslide}V \\
   \bullet @>>{\Rlh}> \bullet
   \end{CD}
   \end{equation*}
\end{prop}
\end{remark} 

%%%%%%%%%%%%%%%%%%%%%%%%%%%%%%%%%%%%%%%%%%%%%%%%%%%%%%%%%%%%%%%%%%%%%%%%%%%%%%%%%%%%%%%%%%%%%%%%%%%%%%%%%%%%%%%%%%%%%%%%%
%
%%%%%%%%%%%%%%%%%%%%%%%%%%%%%%%%%%%%%%%%%%%%%%%%%%%%%%%%%%%%%%%%%%%%%%%%%%%%%%%%%%%%%%%%%%%%%%%%%%%%%%%%%%%%%%%%%%%%%%%%%
\section{Proof of Theorem~\ref{thm: Rslide and Rlb}} \label{pf: Rslide and Rlb}

In this section we give the proof of the following equivalent statement of
Theorem~\ref{thm: Rslide and Rlb}:

Let $\rc \in \Dom(\Rlb^{-1})$ be such that $\Rslide$ is well-defined on $\rc$ and 
$\Rlb^{-1}(\rc) \in \Phi(\Dom(\Pslide)\setminus\PBase)$.
Then $\Rslide(\rc) \in \Dom(\Rlb^{-1})$ and $\Rslide$ is well-defined on
$\Rlb^{-1}(\rc)$ 
and the following diagram commutes:
   \begin{equation*}
   \begin{CD}
   \rc@>{\Rlb^{-1}} >> \bullet \\ 
   @V{\Rslide}VV @VV{\Rslide}V \\
   \bullet @>>{\Rlb^{-1}}> \bullet
   \end{CD}
   \end{equation*}

Without loss of generality we shall assume that $\Rslide$ is not the identity map. 

Let us firstly argue that $\Rslide$ is well-defined on $\Rlb^{-1}(\rc)$. Given that 
$\rc \in \Dom(\Rlb^{-1})$, we know that $\rc$ corresponds to the sequence of rectangles 
$((1,1),(r-1,1),\ldots)$.
Moreover, $\rc = \Rlh^{-1}(\Rlh(\rc),t)$ for some $t \ge r$. (Indeed from the condition that
$\Rlb^{-1}(\rc) \in \Phi((\Dom(\Pslide)\setminus\PBase)$ we can deduce a stronger conclusion
$t > r$, but we only need the weaker statement in the our proof. So we actually proved
a stronger result.) Let $(D^{(n+1)},\ldots, D^{(1)})$ be
the $\Rlh^{-1}$-sequence of $\Rlh(\rc)$ with respect to $t$, we know
$D^{(k)}=0$ for $k < r$. By the definition of $\Rlh^{-1}$,  $\rc^{(k)}$ is obtained
from $\Rlh(\rc)^{(k)}$ by adding a singular string of size 1 for each $k < r$.  

Let $(I^{(n+1)},\ldots, I^{(1)})$ be the $\Rslide$-sequence of $\rc$. We observe that 
$I^{(k)}>1$ for each $k<r$. To see this, let $j$ be the least index such that 
$D^{(j)} > 0$ (the existence of such a $j$ follows from Lemma~\ref{lemma: below} below).
Now $\rc \in \Dom(\Rlb^{-1})$ implies that $j \ge r$. 
Note that all strings in $\rc^{(j)}$ of length $\le D^{(j)}$ are non-singular, 
it follows that the smallest singular string of $\rc^{(j)}$ is at least of length 
$D^{(j)}+1 > 1$. Since $(I^{(n+1)},\ldots, I^{(1)})$ increases, we obtain $I^{(k)}>1$ 
for each $k<r$.

Now $\Rlb^{-1}$ acts on $\rc$ by removing a singular part of size 1 from $\rc^{(k)}$ for
each $k < r$, and leaving the label and colabel of the remaining parts unchanged. Then by
the result from the previous paragraph, the $\Rslide$-sequence of $\Rlb^{-1}(\rc)$ is exactly 
the $\Rslide$-sequence of $\rc$. This shows that $\Rlb^{-1}(\rc) \in \Dom(\Rslide)$.

We extract from the above arguments the following fact, which will be referred to in the future
sections: 
\begin{lemma} \label{lemma: I unchanged under Rlb inverse}
For any $k \in [n+1]$, $\Rlb^{-1}(I)^{(k)} = I^{(k)}$. 
\end{lemma}
 
Let us secondly argue that $\Rslide(\rc) \in \Dom(\Rlb^{-1})$.
We note that $\rc$ and $\Rslide(\rc)$ correspond to the same sequence of rectangles 
$((1,1),(r-1,1),\ldots)$.
$\rc \in \Dom(\Rlb^{-1})$ means that for each $k<r$ $\rc^{(k)}$ has a singular part of size 1,
The arguments from the above paragraph show that $\Rslide$ does not touch these parts, 
thus for each $k<r$, $\Rslide(\rc)^{(k)}$ has a singular part of size 1. 
Therefore $\Rslide(\rc) \in \Dom(\Rlb^{-1})$.

The above arguments also clearly show that $\Rslide(\Rlb^{-1}(\rc)) = \Rlb^{-1}(\Rslide(\rc))$. 

Now the only thing left is the following lemma:
\begin{lemma} \label{lemma: below}
Let $\rc$, $\Rlh(\rc)$ and $D_t=(D^{(n+1)},\ldots, D^{(1)})$ be given as above. Then there 
exists a least index $j$ such that $D^{(j)}>0$.
\end{lemma} 
\begin{proof}
By the definition of $\Rlh^{-1}$ (see Definition~\ref{def: Rl_inverse}), the above statement
is clearly true for $t < n+2$ since $D^{(k)}=\infty$ for $k\ge t$.

In the case that $t=n+2$, our assumption that $\Rlb^{-1}(\rc) \not\in \Phi(\PBase)$ implies that
$\rc \not \in \Phi(\PBase)$. This then implies $\Rlh(\rc)^{(n+1)} \not = \emptyset$. 
By Proposition~\ref{prop: Rslide and Rlh}, $\Rslide$ is well-defined on $\Rlh(\rc)$, in
particular this means that $\Rlh(\rc)^{(n+1)}$ contains a singular string. 
Thus $D^{(n+1)}>0$.
\end{proof}
%%%%%%%%%%%%%%%%%%%%%%%%%%%%%%%%%%%%%%%%%%%%%%%%%%%%%%%%%%%%%%%%%%%%%%%%%%%%%%%%%%%%%%%%%%%%%%%%%%%%%%%%%%%%%%%%%%%%%%%%%
%
%%%%%%%%%%%%%%%%%%%%%%%%%%%%%%%%%%%%%%%%%%%%%%%%%%%%%%%%%%%%%%%%%%%%%%%%%%%%%%%%%%%%%%%%%%%%%%%%%%%%%%%%%%%%%%%%%%%%%%%%%
\section{Proof of Theorem~\ref{thm: Rslide and Rls}} \label{pf: Rslide and Rls}

In this section we give the proof of the following equivalent statement of
Theorem~\ref{thm: Rslide and Rlb}:

Let $\rc \in \Dom(\Rls^{-1})$ be such that $\Rslide$ is well-defined on $\rc$ and
$\Rls^{-1}(\rc) \in \Phi(\Dom(\Pslide)\setminus\PBase)$.
Then $\Rslide(\rc) \in \Dom(\Rls^{-1})$ and $\Rslide$ is well-defined on $\Rls^{-1}(\rc)$,
and the following diagram commutes:
   \begin{equation*}
   \begin{CD}
   \rc@>{\Rls^{-1}} >> \bullet \\ 
   @V{\Rslide}VV @VV{\Rslide}V \\
   \bullet @>>{\Rls^{-1}}> \bullet
   \end{CD}
   \end{equation*}

Let $p=\Phi^{-1}(\rc)$. Then by the condition that $\rc \in \Dom(\Rls^{-1})$,
we have 
\begin{equation*}
   \Yboxdim25pt
    p=S \otimes T \otimes q = \young(\bONE,\vdots,\bR) 
      \otimes \young(\Tonek\cdots\Toneone,\vdots\vdots\vdots,\Trk\cdots\Trone) 
      \otimes q \; ,
\end{equation*}
where $S$ is a single column tableau of height $r$, 
and $T$ is a tableau of shape $(r,s)$ with $s\ge 0$.
By the condition $\Rls^{-1}(\rc) \not\in \PBase$, we have $n+2$ appears in $q$ and
$1$ does not appear in $S$ nor $T$, in particular $\bONE>1$ and $\bR>r$. 

We shall induct on $s$, the number of columns of $T$. To facilitate the induction, let us
denote $\rc_s = \Phi(p)$ where $p$ is defined as above, that is, the subscript $s$ of
$\rc_s$ indicates the width of tableau $T$. 

Let $(I_s^{(n+1)},\ldots,I_s^{(1)})$ be the $\Rslide$-sequence of $\rc_s$. 

The hypothesis we want to carry across inductive steps is the logical  
disjunction of the following two sufficient conditions for the commutativity of above 
diagram:
\begin{hypothesis}[Simplified version] \label{hypo: primitive} 
For each $s\in \mathbb{Z}_{\ge 0}$, $\rc_s$ satisfies one of the following two conditions:
\begin{enumerate}
  \item[A.] $I_s^{(r+1)}>s$;  
  \item[B.] $I_s^{(r+1)} \le s$ but the colabel of any part of $\rc_s^{(r)}$ with sizes in 
        $[I_s^{(r+1)},s]$ is $\ge 2$. 
\end{enumerate} 
\end{hypothesis}

To see that the first condition is sufficient, we recall that $\Rls^{-1}(\rc_s)$ decreases the 
colabels for all parts in $\rc_s^{(r)}$ of size $\le s$ by $1$. 
Now if $I_s^{(r+1)}>s$ then $\Rls^{-1}$ will not affect the choice of $I_s^{(r)}$. In this
case it is also easy to see that $\Rslide$ will not affect the action of $\Rls^{-1}$.

To see that the second condition is sufficient, we notice that if $I_s^{(r+1)} \le s$  
and the colabel of any part of $\rc_s^{(r)}$ with sizes in $[I_s^{(r+1)},s]$ is $\ge 2$, 
then when $\Rls^{-1}$ decreases the colabels for parts in $\rc_s^{(r)}$ of size $\le s$ by $1$,
the parts with size in the interval $[I_s^{(r+1)},s]$ will still be non-singular. Thus 
$\Rls^{-1}$ will not affect the choice of $I_s^{(r)}$. In this case, this colabel 
condition also prevents $\Rslide$ from affecting $\Rls^{-1}$: $\Rslide$ decreases the colabels
of parts in $\rc_s^{(r)}$ of size in the interval $[I_s^{(r+1)},I_s^{(r)})$ by $1$,
thus the parts in $\Rslide(\rc_s)^{(r)}$ of size in $\le s$ will still be non-singular. 

Conditions A and B can be viewed as two possible states of $\rc_s$. We will show below in a more
precise setting that $\rc_s$ either stays in its current state or transits from A 
to B (but never comes back). Intuitively, we can imagine $s$ being time, and when $s$ is small 
($s$ starts from 0), $\rc_s$ starts in state A. As time goes by, we start merging 
columns of height $r$ on the left of $p$, corresponding on the $\RC$ side to $\Rls^{-1}(\rc_s)$,
the condition that after each merging we get a valid tableau (weakly increasing along rows)
corresponds on the $\RC$ side to the fact that $\Rlh^{-1}$ and $\Rlb^{-1}$ affect  
$\rc_s^{(k)}$ for $k>r$ ``less and less''. 
As time pass by, $I_s^{(r+1)}$ will possibly stabilize, and $s$ will possibly catch up and pass 
$I_s^{(r+1)}$, and from this time on $\rc_s$ will get stuck in the state B forever.
  
To make all above statements precise, we first need to set up some notation.

Let $\rc_{j,s}$ for $j\le r$ be the image under $\Phi$ of the following path element:
\begin{equation*} 
   \Yboxdim25pt
S \otimes T \otimes q = 
\young(\bONE,\vdots,\bJ)
       \otimes
\young(\Tonek\cdots\Toneone,%
       \vdots\vdots\vdots,%
       \Trk\cdots\Trone)  \; 
       \otimes
q.
\end{equation*}

Define $\rc_{0,s}=\Phi(T\otimes q)=\Rls^{-1}(\rc_{s-1})$. In particular, $\rc_{0,0}=\Phi(q)$. 
With this notation, our previously defined $\rc_s$ is denoted by $\rc_{r,s}$. 
 
Let $D_{T_{j,s}}=(D_{T_{j,s}}^{(n+1)},\ldots,D_{T_{j,s}}^{(1)})$ be the $\Rlh^{-1}$-sequence 
with respect to $T_{j,s}$. To avoid awkward subindices, we denote $D_{j,s}=D_{T_{j,s}}$.

Let $I_{j,s}=(I_{j,s}^{(n+1)},\ldots,I_{j,s}^{(1)})$ be the $\Rslide$-sequence of $\rc_{j,s}$.
Our previously defined $I_s$, in this notation, is denoted by $I_{r,s}$.
We denote by $I_{0,s}=(I_{0,s}^{(n+1)},\ldots,I_{0,s}^{(1)})$ the $\Rslide$-sequence 
of $\rc_{0,s}=\Rls^{-1}(\rc_{r,s-1})$.

Equipped with this notation, we can give a precise description on how $\Rlh^{-1}$-sequences
interact with $\Rslide$-sequences.

For a fixed $s$, for each $k \in [n+1]$, let us consider the sequence 
$(I_{j,s}^{(k)})_{j\in [r]}$. By Lemma~\ref{lemma: I increases under Rlh inverse},  
we have:
\begin{lemma} \label{lemma: I(s,j) weakly increasing wrt j}
$(I_{j,s}^{(k)})_{j\in [r]}$ is weakly increasing with respect to $j$. More precisely, for any
$s$ and any $k$, we have $I_{j,s}^{(k)} \le I_{j+1,s}^{(k)}$ for $j=0,\ldots,r-1$.
\end{lemma}

Set $A_{j,s}=\{k \mid I_{j,s}^{(k)}>I_{0,s}^{(k)}\}$.
Then by Lemmas~\ref{lemma: I increases under Rlh inverse} 
and~\ref{lemma: same column} we have 
\begin{lemma} \label{lemma: cross}
$A_{j,s}=[1,a_{j,s}]$ for some $1 \le a_{j,s} \le n+1$, and 
$a_{1,s} \le \cdots \le a_{r,s}$. Moreover, 
$A_{j,s}\setminus A_{j-1,s}=\{k \mid I_{j,s}^{(k)}>I_{j-1,s}^{(k)}\}$.
\end{lemma}

An equivalent statement of the above lemma is the following: 
\begin{lemma}
$(I_{j,s}^{(k)})_{j\in [r]}$ satisfies $I_{j,s}^{(k)}=I_{j+1,s}^{(k)}$ everywhere except possibly 
one position $J_{k,s}$ where the inequality is strict $I_{J_{k,s},s}^{(k)}<I_{J_{k,s}+1,s}^{(k)}$.
$J_{k,s}$ is weakly increasing with respect to $k$ and 
$A_{j,s}\setminus A_{j-1,s}=\{k \mid J_{k,s}=j\}$.
\end{lemma}

If we specialize Lemma~\ref{lemma: I(s,j) weakly increasing wrt j} 
to the case $k=r+1$, then $I_{j,s}^{(r+1)}$ is weakly increasing
with respect to $j$. Now we note that $\Rls^{-1}$ does not affect $\rc_s^{(k)}$ for $k \ge r+1$.
Thus $I_{0,s+1}^{(r+1)}=I_{r,s}^{(r+1)}$. Now by the fact that $I_{j,s+1}^{(r+1)}$ is weakly
increasing with respect to $j$, we conclude the following lemma:

\begin{lemma} \label{lemma: I increases at r+1}  
$I_s^{(r+1)}$ is weakly increasing with respect to $s$.
\end{lemma}  

\begin{remark}
The above argument shows that for any fixed $k \ge r+1$, $I_{j,s}^{(k)}$ is weakly increasing
with respect to $(s,j)$ in dictionary order. Indeed this is also the case for any
$k < r+1$, but this fact is a consequence of Theorem~\ref{thm: Rslide and Rls} which we
are going to prove.  
\end{remark}

Now we can finally make precise the statement we want to prove.
\begin{hypothesis}[Precise version] \label{hypo: precise} 
For each $s$, $\rc_s$ satisfies:
\begin{enumerate}
  \item[IH1.] $I_s^{(r)} > s+1$.
  \item[IH2.] $I_{s}^{(r+1)} = I_{s-1}^{(r+1)}$ when $D_{s-1}^{(r)} \ge I_{s-1}^{(r+1)}$.
  \item[IH3.] $\rc_s$ is in one of the following states:
     \begin{itemize}
        \item[SA1 ] $I_s^{(r+1)}>s$, and $D_s^{(r)}<I_s^{(r+1)}$. 
        \item[SA2 ] $I_s^{(r+1)}>s$, and $D_s^{(r)} \ge I_s^{(r+1)}$.  
        \item[SB ] $I_s^{(r+1)} \le s$, and the colabel of any part of $\rc_s^{(r)}$ with size in 
        $[I_s^{(r+1)},s]$ is $\ge 2$, and the colabel of any part with size $s+1$ is $\ge 1$. 
     \end{itemize} 
\end{enumerate} 
\end{hypothesis}

We remark here that [IH3] is the core of the hypothesis, where
SA1 and SA2 combined correspond to the condition A in Hypothesis~\ref{hypo: primitive},
and SB corresponds to the condition SB in Hypothesis~\ref{hypo: primitive}. As we have argued 
before, being in one of these states implies the commutativity of $\Rslide$ and $\Rls^{-1}$. 

By Lemma~\ref{lemma: I increases at r+1}, $I_{s+1}^{(r+1)} \ge I_{s}^{(r+1)}$.
[IH2] then states that $I_s^{(r+1)}$ stabilizes (stops increasing) once $D_s^{(r)} \ge I_s^{(r+1)}$. 

\subsection{Base case:}
We have $s=0$ as our base case.

Let us first verify [IH1], that is, $I_0^{(r)}=I_{r,0}^{(r)}>1$. Suppose not, 
then it must be that $I_{r,0}^{(r)}=1$, and thus $I_{j,0}^{(r)}=1$ for all $j=0,\ldots,r$ 
by the fact that $I_{j,0}^{(r)}$ is weakly increasing with respect to $j$ 
(see Lemma~\ref{lemma: I(s,j) weakly increasing wrt j}). But this is impossible 
since $\bR > r$ and thus $D_{r,0}$ and $I_{r-1,0}$ cross at position $\ge r$ 
(by Lemma~\ref{lemma: cross}, this implies $I_{r,0}^{(r)}>I_{0,0}^{(r)}$, a contradiction).

[IH2] is an empty statement for base case.

Finally, it is clear that $\rc_0$ is either in SA1 or SA2, thus verifies [IH3]. 

\subsection{Induction:} Assume that $\rc_s$ satisfies Hypothesis~\ref{hypo: precise}.
We consider $\rc_{s+1}$.

Let us first verify [IH1], that is, $I_{s+1}^{(r)}>s+2$. 
By [IH1] on $\rc_s$, $I_s^{(r)}>s+1$, so we just need to show $I_{s+1}^{(r)}>I_s^{(r)}$. 
Suppose not, so that $I_{s+1}^{(r)} \le I_s^{(r)}$. By [IH3] on $\rc_s$, $\Rls^{-1}$ does not 
interfere with the $\Rslide$-sequence, so $I_{0,s+1}^{(r)} = I_s^{(r)}$.  
Now by Lemma~\ref{lemma: I(s,j) weakly increasing wrt j} $I_{s+1}^{(r)} \ge I_{0,s+1}^{(r)}$, 
thus $I_{s+1}^{(r)} \ge I_s^{(r)}$. Therefore we it must be the case $I_{s+1}^{(r)} = I_s^{(r)}$. 
This then implies that $I_{j,s+1}^{(r)}=I_s^{(r)}$ for all $j \in [r]$. However, this is 
impossible since $\Lrone = \bR > r$ and thus $D_{r,s+1}$ and $I_{r-1,s+1}$ cross at position 
$\ge r$ (by Lemma~\ref{lemma: cross}, this implies $I_{s+1}^{(r)}>I_s^{(r)}$, a contradiction).

The verification of [IH2] and [IH3] depends on the state of $\rc_s$:

\subsubsection{$\rc_s$ is in state SA1}
By assumption, $D_s^{(r)} < I_s^{(r+1)}$, so [IH2] is vacuously true in this case.
 
We verify [IH3] by showing that $\rc_{s+1}$ is either in SA1 or SA2.
For this we just need to argue that $I_{s+1}^{(r+1)} > s+1$. 
Knowing that $I_{s}^{(r+1)} > s$ (since $\rc_s$ in SA1) and $I_{s+1}^{(r+1)} \ge I_{s}^{(r+1)}$ 
(Lemma~\ref{lemma: I increases at r+1}), it is clear that we just need to show the 
impossibility of $I_{s+1}^{(r+1)} = s+1$. Suppose otherwise. Then it must be the case that 
$I_{s+1}^{(r+1)} = I_s^{(r+1)} = s+1$. This implies that $D_s^{(r)}=s$ (it is clear that
$D_s^{(r)} \ge s$ by the fact that $\rc_s \in \Dom(\Rls^{-1})$, and it is also clear
that $D_s^{(r)} \le s$ by the fact that $D_s^{(r)}<I_s^{(r+1)}$).  
This then implies that $I_s^{(r)}=s+1$ (since $\rc_s$ has a singular string of size 
$D_s^{(r)}+1=s+1$, and $I_s^{(r+1)}=s+1$). But this is a contradiction to [IH1] on $\rc_s$. 

\subsubsection{$\rc_s$ is in state SA2}
Let us verify [IH2] on $\rc_{s+1}$. Given $\rc_s$ in state SA2, we have 
$D_s^{(r)} \ge I_s^{(r+1)}$, so we need to show $I_{s+1}^{(r+1)} = I_s^{(r+1)}$. 
We have argued before that $I_{0,s+1}^{(r+1)}=I_s^{(r+1)}$, so it suffices for us to show
that $I_{r,s+1}^{(r+1)}=I_{0,s+1}^{(r+1)}$. Using Lemmas~\ref{lemma: same column}, 
\ref{lemma: same row} and~\ref{lemma: cross}, it can be inductively shown that 
$D_s^{(r)} \ge I_s^{(r+1)}$ implies $D_{j,s+1}^{(r)}>I_{j-1,s+1}^{(r)}$ and 
$I_{j,s+1}^{(r)}=I_{0,s+1}^{(r)}$ for all $j \in [r]$. 

We verify [IH3] by showing that $\rc_{s+1}$ is either in state SA2 or SB.
Our first observation is that by [IH2] on $\rc_{s+1}$, which we have just verified above, 
$I_{s+1}^{(r+1)}=I_s^{(r)}\le D_s^{(r)}< D_s^{(r+1)}$ (the last $<$ follows from 
Lemma~\ref{lemma: same row}). Thus $\rc_{s+1}$ cannot be in state SA1. 

Next, $\rc_s$ is in state SA2 means that $D_s^{(r)} \ge I_s^{(r+1)} > s$. 
There are two possibilities: 
either $I_s^{(r+1)} > s+1$ or $I_s^{(r+1)} = s+1$. If $I_s^{(r+1)}>s+1$, 
then by [IH2] on $\rc_{s+1}$, we know $I_{s+1}^{(r+1)}=I_s^{(r+1)}>s+1$.
Hence $D_{s+1}^{(r)} \ge I_{s+1}^{(r+1)} > s+1$, so $\rc_{s+1}$ is in SA2. 
If $I_s^{(r+1)} = s+1$, then we want to show that $\rc_{s+1}$ is in state SB. For this
we need to show that in $\rc_{s+1}^{(r)}$ the colabel of any part of size $s+1$ is 
$\ge 2$ and the colabel of any part of size $s+2$ is $\ge 1$.
To see this, we note that in $\rc_s^{(r)}$ the colabel of any part of size $s+1$ is $\ge 1$ 
(otherwise $I_s^{(r)}=s+1$, contradicting [IH1] on $\rc_s$). 
Since $\Rls^{-1}$ does not change the colabel of 
any part in $\rc_s^{(r)}$ of size $s+1$, in $\Rls^{-1}(\rc_s)^{(r)}$ the colabel of any part 
of size $s+1$ is $\ge 1$. Moreover, by Lemma~\ref{lemma: same column}, 
$D_{s+1}^{(r)} > s+1$ implies $D_{k,s+1}^{(r)}> s+1$ for all $k\in [r]$, thus the colabel of any 
part in $\rc_s^{(r)}$ of size $s+1$ is weakly increasing along this sequence of applications of 
$\Rlb^{-1}\circ\Rlh^{-1}$. Furthermore, $D_{r,s+1}^{(r)}> s+1$ implies  that in 
$\rc_{s+1}^{(r)}$ the colabel of any part of size $\le s+2$ increases by $1$. Thus any part 
of size $s+1$ must be of colabel $\ge 2$, and any part of size $s+2$ must be of 
colabel $\ge 1$.
 
\subsubsection{$\rc_s$ is in state SB}
This implies $I_s^{(r+1)} \le s$, which further implies $D_s^{(r)} \ge I_s^{(r+1)}$, 
so [IH2] on $\rc_{s+1}$ is not vacuous in this case. 
The verification is exactly the same as for the case that $\rc_s$ is in state SA2.

Let us verify [IH3] by showing that $\rc_{s+1}$ must be in SB. 
Since $I_{s+1}^{(r+1)}=I_s^{(r+1)}\le s < s+1$, $\rc_{s+1}$ can not be in state SA1 nor SA2. 
Moreover, we note that $\Rls^{-1}$ decreases colabels of all parts in $\rc_s^{(r)}$ of 
size $\le s$, so the colabels of all parts in $\Rls^{-1}(\rc_s)^{(r)}$ of size $\le s+1$ 
are $\ge 1$. Then by Lemma~\ref{lemma: same column}, $D_{s+1}^{(r)} > s+1$ implies that  
$D_{k,s+1}^{(r)}> s+1$ for all $k\in [r]$. Hence the colabel of any 
part in $\rc_s^{(r)}$ of size $\le s+1$ is weakly increasing along this sequence of 
applications of $\Rlb^{-1}\circ\Rlh^{-1}$. Furthermore, $D_{r,s+1}^{(r)}> s+1$ implies  that in 
$\rc_{s+1}^{(r)}$ the colabel of any part of size $\le s+2$ increases by $1$, thus any part 
of size $s+1$ must be of colabel $\ge 2$, and any part of size $s+2$ must be of 
colabel $\ge 1$.

This finish the induction.

%%%%%%%%%%%%%%%%%%%%%%%%%%%%%%%%%%%%%%%%%%%%%%%%%%%%%%%%%%%%%%%%%%%%%%%%%%%%%%%%%%%%%%%%%%%%%%%%%%%%%%%%%%%%%%%%%%%%%%%%%
% 
%%%%%%%%%%%%%%%%%%%%%%%%%%%%%%%%%%%%%%%%%%%%%%%%%%%%%%%%%%%%%%%%%%%%%%%%%%%%%%%%%%%%%%%%%%%%%%%%%%%%%%%%%%%%%%%%%%%%%%%%%
\section{Proof of Theorem \ref{thm: Pslide and Prs}} \label{pf: Pslide and Prs}
   Let $p \in \PBase \setminus \PBBaseA$. Then Theorem~\ref{thm: Pslide and Prs} states that 
   $\Prs(p) \in \PBase \setminus \PBBaseA$ and $\Pslide \circ \Prs(p) = \Prs \circ \Pslide(p)$.

The main idea of the proof is to observe that both $\Pslide$ and $\Prs$ are ``local'' operations,
and their actions are ``far away'' from each other. Thus one operation will not interfere 
the other operation. The notion of ``local'' and ``far away'' are made precise below.

By assumption, $p=T_0 \otimes \cdots \otimes T_K = T_0 \otimes M =
N \otimes T_K$ with $n+2$ appearing only in $T_0$. 
By definition, $\Pslide$ acts on $T_0$ and the changes made to $T_0$ do not depend 
on $M$. In another words, $\Pslide(T_0 \otimes M) = \Pslide(T_0) \otimes M$. 
In this sense, $\Pslide$ acts ``locally'' on the left.

We claim that $\Prs$ changes only $T_K$, and the changes made to $T_K$ do not depend on
$N$. In another word, $\Prs(N \otimes T_K) = N \otimes \Prs(T_K)$. Thus $\Prs$ acts ``locally''
on the right.  
Given that $p \not\in \PBBaseA$, $T_0$ and $T_K$ are distinct. Hence 
$\Pslide$ and $\Prs$ act on distinct tensor factors and therefore ``far away'' from
each other.

Therefore, we are left to show $\Prs(N \otimes T_K) = N \otimes \Prs(T_K)$. This is clearly 
implied by the following general statement: 
\begin{lemma} \label{lemma: Prs and lh/lb/ls}
Let $p = T_1 \otimes \cdots \otimes T_K \in \Path$ where $K > 1$.
Then $\Prs$ commutes with $\Plh$ (or $\Plb$ or $\Pls$) for $p \in \Dom(\Plh)$ 
(or $p \in \Dom(\Plb)$ or $p \in \Dom(\Pls)$, respectively):
   \begin{equation*}
   \begin{CD}
   p@>{\Plh/\Plb/\Pls}>> \bullet \\ 
   @V{\Prs}VV @VV{\Prs}V \\
   \bullet @>>{\Plh/\Plb/\Pls}> \bullet
   \end{CD}
   \end{equation*}
\end{lemma}

To prove Lemma~\ref{lemma: Prs and lh/lb/ls}, we recall that $\Prs$ is defined to be 
$\Phi^{-1}\circ\Rrs\circ\Phi$, so it is natural to consider the following counterpart statement 
on the $\RC$ side: 
\begin{lemma} \label{lemma: Rrs and lh/lb/ls}
Let $\rc \in \RC(B)$, where $B=((r_1,s_1), \ldots, (r_K,s_K))$ where $K > 1$.
Then $\Rrs$ commutes with $\Rlh$ (or $\Rlb$ or $\Rls$) for $\rc \in \Dom(\Rlh)$ 
(or $\rc \in \Dom(\Rlb)$ or $\rc \in \Dom(\Rls)$, respectively). 
   \begin{equation*}
   \begin{CD}
   \rc@>{\Rlh/\Rlb/\Rls}>> \bullet \\ 
   @V{\Rrs}VV @VV{\Rrs}V \\
   \bullet @>>{\Rlh/\Rlb/\Rls}> \bullet
   \end{CD}
   \end{equation*}
\end{lemma}
\begin{proof}
Recall the action of $\Rrs$ on $\rc$: If $\rc \in \RC(B)$, then 
$\Rrs(\rc) \in \RC(\rs(B))$ is obtained by increasing the labels by $1$ for all 
parts in $\rc^{(r_K)}$ of size less than $s_K$. In particular, $\Rrs$ leaves the 
colabels of all parts unchanged. 

By definition, the action of $\Rlh$ depends only on the colabels of parts, 
and $\Rlh$ preserves labels for all unchanged parts, thus commutes with $\Rrs$.

$\Rlb$ adds a singular string of length 1 to each $\rc^{(a)}$ for $1 \le a < r_1$ 
and preserves both labels and colabels for all other parts, while 
$\Rrs$ preserves colabels and thus the singularity of all parts. So $\Rlb$ and $\Rrs$ commute.

The action of $\Rls$ splits one column from left of the rectangle $(r_1,s_1)$, and
increases the colabel of any part of size $< s_1$ in $\rc^{(r_1)}$ by $1$.
The action of $\Rrs$ splits one column from right on the rectangle $(r_K,s_K)$, and
increases the label of any part of size $< s_K$ in $\rc^{(r_K)}$ by $1$.
Clearly they commute.
\end{proof}

To see that Lemma \ref{lemma: Rrs and lh/lb/ls} implies Lemma \ref{lemma: Prs and lh/lb/ls},
we consider now the following diagram:
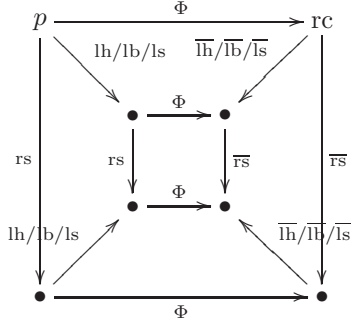
\begin{equation*}
\xymatrix{
 {p} \ar[rrr]^{\Phi} \ar[ddd]_{\Prs} \ar[dr] ^{\Plh/\Plb/\Pls}& & &
        {\rc} \ar[ddd]^{\Rrs} \ar[dl]_{\Rlh/\Rlb/\Rls} \\
 & {\bullet} \ar[r]^{\Phi} \ar[d]_{\Prs} & {\bullet} \ar[d] ^{\Rrs}& \\
 & {\bullet} \ar[r]^{\Phi}  & {\bullet}  & \\
 {\bullet} \ar[rrr]_{\Phi} \ar[ur] ^{\Plh/\Plb/\Pls}& & & {\bullet} \ar[ul]_{\Rlh/\Rlb/\Rls}
}
\end{equation*}

We observe the following on this diagram:
\begin{enumerate} 
  \item The upper and lower face commutes by the definition of $\Phi$. 
  \item The front and back face commutes by the definition of $\Prs$. 
  \item The right face commutes by Lemma \ref{lemma: Rrs and lh/lb/ls} 
\end{enumerate} 

The above observations imply that the left face commutes which is the statement of 
Lemma \ref{lemma: Prs and lh/lb/ls}. 
This proves the main statement.

%%%%%%%%%%%%%%%%%%%%%%%%%%%%%%%%%%%%%%%%%%%%%%%%%%%%%%%%%%%%%%%%%%%%%%%%%%%%%%%%%%%%%%%%%%%%%%%%%%%%%%%%%%%%%%%%%%%%%%%%%
%
%%%%%%%%%%%%%%%%%%%%%%%%%%%%%%%%%%%%%%%%%%%%%%%%%%%%%%%%%%%%%%%%%%%%%%%%%%%%%%%%%%%%%%%%%%%%%%%%%%%%%%%%%%%%%%%%%%%%%%%%%
\section{Proof of Theorem~\ref{thm: base case 1}} \label{pf: base case 1}
Let $p \in \PBBaseA$, then the following diagram commutes:
\begin{equation*}
\begin{CD}
p@>{\Phi}>> \bullet \\ 
@V{\Pslide}VV @VV{\Rslide}V \\
\bullet @>>{\Phi}> \bullet
\end{CD}
\end{equation*}

Before we start the proof, we first give an alternative description of $\Phi$
for the special case of a single tensor factor.

\subsection{An algorithm for computing $\Phi(p)$ for $p \in \Path(r,c)$}
\begin{definition} \label{def: new Phi}
Given $p\in\Path_n(r,c)$, define $\Psi(p)=(u_1,\ldots,u_n)$, where each 
$u_k$ is the part of $p$ that is formed by all boxes that
\begin{enumerate} 
  \item are on $j$-th row for $j \le k$;
  \item contain a number $>k$. 
\end{enumerate} 
\end{definition}

\begin{example} \label{example: area}
Let $p \in \Path_7(4,4)$ be $\young(1234,2445,3666,5778)$, then 
\begin{equation*}
\Psi(p)\;=\;
(\young(234),\young(:34,445),\young(::4,445,666),\young(:::5,:666,5778),\young(666,778),\young(778),\young(8)) \; .
\end{equation*}
\end{example}

\begin{remark} \label{rmk: area}
Each $u_k$ in the above list has the shape of a Ferrers diagram of some partition rotated by $180^o$.
This follows from the fact that $p$ is a semi-standard Young tableau. Thus, we can set $\Psit(p)$ 
to be the list of (rotated) partitions associated to $\Psi(p)$. 
It is often useful to think of each $u_k$ literally as the ``area" of $p$ occupied by the boxes
of $u_k$. 
For example, the box $\young(2)$ in $u_1=\young(234)$ is the $(1,2)$-box of $p$; the 
box $\young(3)$ in $u_1=\young(234)$ is the same box containing $3$ in $u_2=\young(:34,445)$ 
since they are both the $(1,3)$-box of $p$. As another example, row $\young(666)$ in 
$u_3$, $u_4$, and $u_5$ literally denotes the $(3,2)$, $(3,3)$ and $(3,4)$-boxes of $p$
in all three cases.

To stress the point we are making, we rewrite Example~\ref{example: area} by highlighting
the boxes of each $u_k$ in $p$ in red:
\begin{equation*}
\Yboxdim10pt
\Psi(p)\;=\;
(\young(1\two\three\four,2445,3666,5778),
 \young(12\three\four,2\four\four\five,3666,5778),
 \young(123\four,2\four\four\five,3\six\six\six,5778),
 \young(1234,244\five,3\six\six\six,\five\seven\seven\eight),
 \young(1234,2445,3\six\six\six,5\seven\seven\eight),
 \young(1234,2445,3666,5\seven\seven\eight),
 \young(1234,2445,3666,577\eight)) \; .
\end{equation*}
\end{remark}

There is an alternative description of $\Psi$, which recursively constructs $\Psi(p)^{(k)}$ 
from $\Psi(p)^{(k+1)}$. Our proof exploits this construction. 
\begin{definition} \label{def: alternative of Psi} \mbox{}
\begin{enumerate}
   \item $\Psi(p)^{(n)}$ is the area of boxes of $p$ that contains $n+1$. The area is clearly
         a horizontal strip.
   \item $\Psi(p)^{(k)}$ is obtained from $\Psi(p)^{(k+1)}$ by adding all boxes of $p$ 
         that contain $k+1$ (which forms a horizontal strip) and then 
         removing the $(k+1)$-st row of $p$ if $k+1 \le r$. 
\end{enumerate}
\end{definition} 
The equivalence of the above two descriptions is clear.  

We have the following result that relates $\Psit$ and $\Phi$: 
\begin{theorem} \label{thm: Psi=Phi} 
For $p\in \Path(r,c)$, we have $\Psit(p) = \Phi(p)$, where all strings of $\Psit(p)$ are singular.
\end{theorem} 
\begin{proof}
This is proved in Appendix~\ref{pf: Psi=Phi}. 
\end{proof}

\begin{corollary}
For $p\in \Path(r,c)$, the rigged configuration $\Phi(p)$ has only singular strings.
\end{corollary}

From now on, we identify $\Psit(p)$ with a rigged configuration as
described in Theorem~\ref{thm: Psi=Phi}. 

\subsection{Jeu-de-taquin on $\Psi(p)$}
By Remark~\ref{rmk: area}, each element $u_k$ in $\Psi(p)$ can be viewed as a collection
of boxes in $p$. Thus jeu-de-taquin on $p$ is directly reflected on $\Psi(p)$. 

\begin{definition} \label{def: Aslide}
For $p \in \Path_{n+1}(r,c)$, define $\Aslide(\Psi(p))=\Psi(\rho(p))$.
\end{definition}

\begin{example} \label{example.rhobb}
\begin{multline*}
\Aslide(\Psi(p))\;= \\
(\young(11\three\four,2245,3466,5677),
 \young(11\three\four,22\four\five,3466,5677),
 \young(113\four,22\four\five,3\four\six\six,5677),
 \young(1134,224\five,34\six\six,\five\six\seven\seven),
 \young(1134,2245,34\six\six,5\six\seven\seven),
 \young(1134,2245,3466,56\seven\seven),
 \young(1134,2245,3466,5677)) \; .
\end{multline*}
\end{example}

Let $SR$ be the sliding route of $p \in \PBBaseA$ under $\Pslide$ 
(see Definition~\ref{def: Pslide}). We have:
\begin{lemma}
$SR$ intersects $\Psi(p)^{(k)}$ for each $k$. If $(i,j)$ is the last box in $SR$ (recall 
$SR$ is a sequence of boxes of $p$) that 
is also in $\Psi(p)^{(k)}$ then $(i,j)$ is an upper left corner of $\Psi(p)^{(k)}$. 
\end{lemma}
\begin{proof}
We shall induct on the recursive definition of $\Psi$ 
(see Definition~\ref{def: alternative of Psi}) with the following hypothesis:
\begin{hypothesis}
For each $k \in [n+1]$, 
\begin{enumerate}
   \item [IH1.] $SR$ intersects $\Psi(p)^{(k)}$. If $(i_k,j_k)$ is the last box in $SR$ that 
         is also in $\Psi(p)^{(k)}$ then $(i_k,j_k)$ is an upper left corner of $\Psi(p)^{(k)}$. 
   \item [IH2.] $i_k \le k$. 
\end{enumerate}
\end{hypothesis}

In the base case we consider $\Psi(p)^{(n+1)}$. Given that $p \in \PBBaseA$, there is a 
horizontal strip of ($n+2$)'s in $p$, which forms $\Psi(p)^{(n+1)}$. By the definition of 
sliding route, an initial segment of $SR$ overlaps with this horizontal strip, 
thus $i_{n+1}=r$. It is clear that IH1 and IH2 hold. 

Now assume that IH1 and IH2 hold for $\Psi(p)^{(k+1)}$.  
From $\Psi(p)^{(k+1)}$ to $\Psi(p)^{(k)}$ a horizontal strip of boxes that contains $k+1$ is 
added. We distinguish two cases:
\begin{enumerate}  
  \item There is a box containing $k+1$ above box $(i_{k+1},j_{k+1})$ of $p$.
  \item There is no box containing $k+1$ above box $(i_{k+1},j_{k+1})$ of $p$.
\end{enumerate}  

In the first case, let $(i_k, j_k)=(i_{k+1}-1, j_k)$. Then IH1 and IH2 hold for $\Psi(p)^{(k)}$.  

In the second case, let $(i_k, j_k)=(i_{k+1}, j_k-l)$ where $l$ is the length of the 
horizontal strip of $(k+1)$'s in row $i_k$ (possibly $l=0$). To check IH2, we assume the 
contrary that $i_k=i_{k+1}>k$. 
We also have $p_{i_k-1,j_k}\le k$. Thus $p_{i_k-1,j_k} \le i_k-1$. 
This is impossible for $p \in \PBBaseA \subset \Dom(\Pslide)$.
IH2 then implies that the $i_k$-th row will not be removed while constructing 
$\Psi(p)^{(k)}$ from $\Psi(p)^{(k+1)}$, hence IH1 holds. 
\end{proof}

By the above lemma we see that $\Aslide$ removes a box from some part of size $b_k$ in
$\Psi(p)^{(k)}$ for each $k$. We call $(b_{n+1}, \ldots, b_1)$ 
the \textbf{$\Aslide$-sequence} of $\Psi(p)$ (or equivalently $\Phi(p)$).

\begin{example}
Continuing Example~\ref{example.rhobb} we have $(b_7,\ldots,b_1) = (1,3,3,3,3,3,3)$.
\end{example}

\subsection{Proof of Theorem~\ref{thm: base case 1}}
As the result of the above two sections, we just need to show $\Rslide(\Psi(p))=\Aslide(\Psi(p))$ 
for $p \in \PBBaseA$. Let $\vec{a}=(a_{n+1},\ldots,a_1)$ be the $\Rslide$-sequence
of $\Phi(p)$ (see Definition~\ref{def: Rslide}).
Let $\vec{b}=(b_{n+1},\ldots,b_1)$ be the $\Aslide$-sequence of $\Psi(p)$ defined above.

It suffices to show that $\vec{a}=\vec{b}$. 
We shall induct on the recursive definition of $\Psi$. 

As base case, $a_{n+1}=b_{n+1}$ since $\Psi(p)^{(n+1)}$ has only one (singular) part. 
Now assume $a_{k+1}=b_{k+1}$. By the definition of the $\Rslide$-sequence, $a_k$ is the size of the 
shortest part (all parts are singular) in $\Psi(p)^{(k)}$ that is no shorter than $a_{k+1}$. 
From $\Psi(p)^{(k+1)}$ to $\Psi(p)^{(k)}$ a horizontal strip of boxes that contains $k+1$ is 
added. There are three possibilities:
\begin{enumerate}  
  \item There is no added box containing $k+1$ adjacent (above or to the left) to the box 
        removed from $\Psi(p)^{(k+1)}$ by $\Aslide$.
  \item There is a box containing $k+1$ added above the box removed from $\Psi(p)^{(k+1)}$ by 
        $\Aslide$. 
  \item There is a box containing $k+1$ added on the left to the box removed from 
        $\Psi(p)^{(k+1)}$ by $\Aslide$. But no box containing $k+1$ added above.
\end{enumerate}  

In the first case, the removed box is already an inside corner of $\Psi(p)^{(k)}$, 
thus by the definition of $\Aslide$, no further Sch{\"u}tzenberger's sliding can be done, 
hence $b_k=b_{k+1}$. It is also clear that the part in $\Psi(p)^{(k)}$ that contains the 
removed box is the shortest part which is no shorter than $a_{k+1}$, thus 
$a_k=a_{k+1}=b_{k+1}=b_k$. 

In the second case, we need one more Sch{\"u}tzenberger's sliding up to get to the 
inside corner when constructing $\Aslide(\Psi(p))^{(k)}$, thus $a_k=a_{k+1}=b_{k+1}=b_k$. 

In the third case, we possibly need several more Sch{\"u}tzenberger's slidings left to 
get to the inside corner when constructing $\Aslide(\Psi(p))^{(k)}$. Then $b_k$ is the 
size of the part that contains the removed box in $\Psi(p)^{(k)}$, but this part is 
clearly also the shortest part which is no shorter than $a_{k+1}$. Therefore, $a_k=b_k>a_{k+1}
=b_{k+1}$.

%%%%%%%%%%%%%%%%%%%%%%%%%%%%%%%%%%%%%%%%%%%%%%%%%%%%%%%%%%%%%%%%%%%%%%%%%%%%%%%%%%%%%%%%%%%%%%%%%%%%%%%%%%%%%%%%%%%%%%%%%
%
%%%%%%%%%%%%%%%%%%%%%%%%%%%%%%%%%%%%%%%%%%%%%%%%%%%%%%%%%%%%%%%%%%%%%%%%%%%%%%%%%%%%%%%%%%%%%%%%%%%%%%%%%%%%%%%%%%%%%%%%%
\section{Proof of Theorem~\ref{thm: base case 2}} \label{pf: base case 2}
In this section we prove that for $p\in \PBBaseB$, we have 
$\Rslide \circ \Phi(p) = \Phi \circ \Pslide(p)$.

Given $p \in \PBBaseB$, we have 
\begin{equation*} 
\Yboxdim25pt
p = S \otimes q = \young(\sONE,\vdots,\sT) \otimes q
\end{equation*} 
where $s_1 > 1$, $s_t=n+2$, and $q \in \Path_{n+1}(B)$ for some $B$, and $n+2$ does not 
appear anywhere in $q$.

For $k \le t$, denote $S_k$ be the single column tableau formed by the first $k$ boxes of $S$. 
That is
\begin{equation*}
\Yboxdim25pt
S_i=\young(\sONE,\sTWO,\vdots,\sI) \;.
\end{equation*}

Then $(\Pslide(S))_{k+1}$ is the single column tableau 
\begin{equation*}
\Yboxdim25pt 
(\Pslide(S))_{k+1}=\young(1,\sONE,\vdots,\sI) \; .
\end{equation*}

Let us first lay out the road map of the proof. We will describe a combinatorial construction 
$\alpha : \RC(((k,1),B),\la) \to \RC(((k+1,1),B),\la+\epsilon_1)$, and inductively argue 
that $\Phi(S_k\otimes q) \stackrel{\alpha}{\mapsto} \Phi((\Pslide(S))_{k+1}\otimes q)$ for all 
$k<t$. In particular, $\Phi(S_{t-1}\otimes q) \stackrel{\alpha}{\mapsto} \Phi(\Pslide(S)\otimes q)$.
Then we argue the commutativity of the following diagram:
   \begin{equation} \label{eq: triangle}
   \xymatrix{
    {\Phi(S_{t-1}\otimes q)} \ar[rr]^{\Rlb^{-1}\circ\Rlh^{-1}(s_t)} \ar[d]^{\alpha} && {\Phi(S\otimes q)} 
                                                                            \ar[d]^{\Rslide}  \\
    {\Phi(\Pslide(S)\otimes q)} \ar[rr]^{=} && {\bullet} 
   } 
   \end{equation}
This then implies $\Phi(\Pslide(S\otimes q))=\Phi(\Pslide(S)\otimes q)=\Rslide(\Phi(S\otimes q))$
which finishes the proof. 

Recall that by Remark~\ref{rmk: colabel}, we can either describe an $\rc \in \RC$ as $(v,J)$ 
in terms of its label or as $(v,cJ)$ in terms of its colabel. In the following definition,
it is more convenient for us to use colabels.

\begin{definition}
Let $k < t$. For $i \in [k]$ let $D_i:=D_{s_i}$ be the $\Rlh^{-1}$-sequence of $\Phi(S_{i-1}\otimes q)$ 
with respect to $s_i$ (see Definition~\ref{def: Rl_inverse}). 
By construction, $D_i^{(i)} \le D_{i-1}^{(i-1)}$. For notational convenience, we take 
$D_0^{(0)}=\infty$. 

Let $E_i$ be the sequence of singular strings in $\Phi(S_i \otimes q)$  that were obtained 
by the action of $\Rlh^{-1}$ that adds a box to each string in $D_i$.    

Let $(v,cJ)=\Phi(S_k\otimes q)$. Then $(\widetilde{v},\widetilde{cJ})=\alpha(v,cJ)$ is defined 
by the following construction (recall Remark~\ref{rmk: succ and >}, when we compare
parts below, we are comparing their length):
\begin{enumerate}  
  \item for $j > k+1$, $(\widetilde{v},\widetilde{cJ})^{(j)}=(v,cJ)^{(j)}$. 
  \item for $j = k+1$, $\widetilde{v}^{(k+1)}=v^{(k+1)}$, $\widetilde{cJ}^{(k+1)}$ is obtained from 
        $cJ^{(k+1)}$ by
        \begin{itemize} 
           \item for strings $s \le D_k^{(k)}$, $\widetilde{cJ}^{(k+1)}(s)=cJ^{(k+1)}(s)+1$;
           \item for strings $s > D_k^{(k)}$, $\widetilde{cJ}^{(k+1)}(s)=cJ^{(k+1)}(s)$. 
        \end{itemize} 
  \item for $j = k$, $\widetilde{v}^{(k)}$ removes one box from the part $E_k^{(k)}$ in $v^{(k)}$.
        $\widetilde{cJ}^{(k)}$ is such that
        \begin{itemize} 
           \item the shortened part has colabel 0;
           \item for strings $s \le D_k^{(k)}$, $\widetilde{cJ}^{(k)}(s)=cJ^{(k)}(s)-1$;
           \item for strings $D_k^{(k)} < s \le D_{k-1}^{(k-1)}$, 
                 $\widetilde{cJ}^{(k)}(s)=cJ^{(k)}(s)+1$; 
           \item for strings $s > D_{k-1}^{(k-1)}$, $\widetilde{cJ}^{(k)}(s)=cJ^{(k)}(s)$. 
        \end{itemize} 
  \item for $1 \le j < k$, $\widetilde{v}^{(j)}$ removes one box from the part 
        $E_j^{(j)}$ in $v^{(j)}$. $\widetilde{cJ}^{(j)}$ is such that
        \begin{itemize} 
           \item the shortened part has colabel 0;
           \item for strings $s \le D_{j+1}^{(j+1)}$, $\widetilde{cJ}^{(j)}(s)=cJ^{(j)}(s)$;
           \item for strings $D_{j+1}^{(j+1)} < s \le |D_{j}^{(j)}|$, 
                 $\widetilde{cJ}^{(j)}(s)=cJ^{(j)}(s)-1$;
           \item for strings $D_{j}^{(j)} < s \le D_{j-1}^{(j-1)}$, 
                 $\widetilde{cJ}^{(j)}(s)=cJ^{(j)}(s)+1$;
           \item for strings $s > D_{j-1}^{(j-1)}$, $\widetilde{cJ}^{(j)}(s)=cJ^{(j)}(s)$. 
        \end{itemize} 
\end{enumerate}  
 
\end{definition}

\begin{lemma} \label{lemma: alpha}
For each $k<t$, $\Phi(S_k\otimes q) \stackrel{\alpha}{\mapsto} \Phi((\Pslide(S))_{k+1}\otimes q)$.
Moreover, $E_k^{(k)}$ is the smallest singular string in $(\Phi(S_k\otimes q))^{(k)}$ of 
size greater than 0; and for $1 \le j \le k-1$, $E_{j}^{(j)}$ is the smallest singular string in 
$(\Phi(S_{j}\otimes q))^{(j)}$ of size greater than $E_{j+1}^{(j+1)}$.
\end{lemma}    
\begin{proof}
Proceed by induction.

In the base case $k=1$, let $(v_1,cJ_1)=\Phi(\young(\sONE)\otimes q)$. By assumption, $\sONE>1$,
thus $D_1^{(1)}<\infty$. 
The definition of $\alpha$ then says $(\widetilde{v_1},\widetilde{cJ_1})=\alpha(v_1,cJ_1)$ is such that 
\begin{enumerate}
   \item $(\widetilde{v_1},\widetilde{cJ_1})^{(j)}=(v_1,cJ_1)^{(j)}$ for $j>2$; 
   \item $\widetilde{v_1}^{(2)} = v_1^{(2)}$, 
         $\widetilde{cJ_1}^{(2)}(s)=cJ_1^{(2)}(s)+1$ for strings $s \le D_1^{(1)}$, 
         $\widetilde{cJ_1}^{(2)}(s)=cJ_1^{(2)}(s)$ for strings $s > D_1^{(1)}$;
   \item $\widetilde{v_1}^{(1)}$ removes one box from the part $E_1^{(1)}$ in $v_1^{(1)}$, 
         $\widetilde{cJ_1}^{(1)}$
         is such that the shortened part has colabel 0, $\widetilde{cJ_1^{(1)}}(s)=cJ_1^{(1)}(s)-1$ for
         strings $s \le D_1^{(1)}$,  $\widetilde{cJ_1^{(1)}}(s)=cJ_1^{(1)}(s)+1$ for
         strings $s > D_1^{(1)}$.  
\end{enumerate}

A direct computation shows that $\alpha(v_1,cJ_1)=\Phi(\young(1,\sONE)\otimes q)$. 
Moreover, the fact that $\widetilde{cJ_1^{(1)}}(s)=cJ_1^{(1)}(s)-1$ for strings $s\le D_1^{(1)}$
implies that $E_1^{(1)}$ is the smallest singular string of size greater than 0 in 
$(\Phi(\young(\sONE) \otimes q))^{(1)}$. This proves the base case.

Now let $(v_k,cJ_k)=\Phi(S_k\otimes q)$ and 
$(\widetilde{v}_k,\widetilde{cJ}_k)=\Phi((\Pslide(S))_{k+1}\otimes q)$,
and suppose that $(v_k,cJ_k)\stackrel{\alpha}{\mapsto} (\widetilde{v_k},\widetilde{cJ_k})$. 
Consider the difference between 
$(v_{k+1},cJ_{k+1})=\Rlb^{-1}\circ\Rlh^{-1}(\Phi(S_k\otimes q),s_{k+1})$
and $(\widetilde{v}_{k+1},\widetilde{cJ}_{k+1})=\Rlb^{-1}\circ\Rlh^{-1}(\Phi((\Pslide(S))_{k+1}
\otimes q),s_{k+1})$.
We will argue that the difference is exactly the effect of $\alpha$.

Let $D_{k+1}$ and $\widetilde{D}_{k+1}$ be the $\Rlh^{-1}$-sequences of $(v_k,cJ_k)$ and 
$(\widetilde{v}_k,\widetilde{cJ}_k)$ with respect to $s_{k+1}>k+1$, respectively. 
$E_{k+1}$ and $\widetilde{E}_{k+1}$ are defined as before corresponding to $D_{k+1}$ and 
$\widetilde{D}_{k+1}$, respectively.  

Consider the difference between $(v_{k+1},cJ_{k+1})^{(a)}$ and 
$(\widetilde{v}_{k+1},\widetilde{cJ}_{k+1})^{(a)}$ for $a>k+2$.
By induction, for $j > k+1$, $(v_k,cJ_k)^{(j)}$ and $(\widetilde{v}_k,\widetilde{cJ}_k)^{(j)}$ 
are the same, thus $D_{k+1}^{(j)}=\widetilde{D}_{k+1}^{(j)}$ for $j=n+1$ down to $k+2$. 
This then implies that for $j > k+2$, 
$(v_{k+1},cJ_{k+1})^{(j)}=(\widetilde{v}_{k+1},\widetilde{cJ}_{k+1})^{(j)}$.

Consider the difference between $(v_{k+1},cJ_{k+1})^{(k+2)}$ and 
$(\widetilde{v}_{k+1},\widetilde{cJ}_{k+1})^{(k+2)}$.
The arguments from the previous paragraph also shows that 
$v_{k+1}^{(k+2)}=\widetilde{v}_{k+1}^{(k+2)}$. 
By induction, $v_k^{(k+1)}=\widetilde{v}_k^{(k+1)}$, 
and for strings $s \le D_{k}^{(k)}$, $\widetilde{cJ}_{k}^{(k+1)}(s)=cJ_{k}^{(k+1)}(s)+1$. 
Then by definition of $\Rlh^{-1}$ we have $D_{k+1}^{(k+1)} \ge \widetilde{D}_{k+1}^{(k+1)}$. 
By the fact that $\Pslide(S)_{k+1}$ is of height $k+1$, we have $\widetilde{D}_{k+1}^{(k+1)}=0$. 
Thus by the definition of $\Rlb^{-1}$, $v_{k+1}^{(k+1)}$ has one box added on the string 
$D_{k+1}^{(k+1)}$ of $v_{k}^{(k+1)}$, while $\widetilde{v}_{k+1}^{(k+1)}=\widetilde{v}_k^{(k+1)}$. 
Therefore, $\widetilde{v}_{k+1}^{(k+1)}$ can be obtained from $v_{k+1}^{(k+1)}$ by removing 
one box from the part $E_{k+1}^{(k+1)}$. All above and the fact that the sequence of rectangles
$((k+2,1),B)$ of $(\widetilde{v}_{k+1},\widetilde{cJ}_{k+1})$ contributes $1$ more to the vacancy number
of strings in $(\widetilde{v}_{k+1},\widetilde{cJ}_{k+1})^{(k+2)}$ than the sequence of rectangles
$((k+1,1),B)$ of $(v_{k+1},cJ_{k+1})$ contributes to the vacancy number of strings in 
$(v_{k+1},cJ_{k+1})^{(k+2)}$ implies that
\begin{itemize} 
   \item for strings $s \le D_{k+1}^{(k+1)}$, $\widetilde{cJ}_{k+1}^{(k+2)}(s)=cJ_{k+1}^{(k+2)}(s)+1$;
   \item for strings $s > D_{k+1}^{(k+1)}$, $\widetilde{cJ}_{k+1}^{(k+2)}(s)=cJ_{k+1}^{(k+2)}(s)$. 
\end{itemize} 

Consider the difference between $(v_{k+1},cJ_{k+1})^{(k+1)}$ and 
$(\widetilde{v}_{k+1},\widetilde{cJ}_{k+1})^{(k+1)}$.
We have just shown in the previous paragraph that $\widetilde{v}_{k+1}^{(k+1)}$ 
is obtained from $v_{k+1}^{(k+1)}$ 
by removing one box from the part $E_{k+1}^{(k+1)}$. By definition of $\Rlb^{-1}$, 
the part with this box removed has colabel 0. By induction, $\widetilde{cJ}_k^{(k+1)}= cJ_k^{(k+1)}+1$
for parts $\le D_k^{(k)}$. Also by induction $\widetilde{v}_{k}^{(k)}$ can be obtained from 
$v_{k}^{(k)}$ by removing one box from the part $E_k^{(k)}$. All above and the fact that the 
sequence of rectangles $((k+2,1),B)$ of $(\widetilde{v}_{k+1},\widetilde{cJ}_{k+1})$ contributes $1$ 
less to the vacancy number of strings in $(\widetilde{v}_{k+1},\widetilde{cJ}_{k+1})^{(k+1)}$ than the 
sequence of rectangles $((k+1,1),B)$ of $(v_{k+1},cJ_{k+1})$ contributes to the vacancy number 
of strings in $(v_{k+1},cJ_{k+1})^{(k+1)}$ implies that
\begin{itemize} 
  \item for strings $s \le D_{k+1}^{(k+1)}$, $\widetilde{cJ}_{k+1}^{(k+1)}(s)=cJ_{k+1}^{(k+1)}(s)-1$;
  \item for strings $D_{k+1}^{(k+1)} < s \le D_k^{(k)}$, 
        $\widetilde{cJ}_{k+1}^{(k+1)}(s)=cJ_{k+1}^{(k+1)}(s)+1$; 
  \item for strings $s > D_k^{(k)}$, $\widetilde{cJ}_{k+1}^{(k+1)}(s)=cJ_{k+1}^{(k+1)}(s)$. 
\end{itemize} 

By the first bullet point above, $E_{k+1}^{(k+1)}$ is the smallest singular string in 
$(v_{k+1},cJ_{k+1})^{(k+1)}$ of size $\ge 0$. 

Consider the difference between $(v_{k+1},cJ_{k+1})^{(k)}$ and 
$(\widetilde{v}_{k+1},\widetilde{cJ}_{k+1})^{(k)}$.
By induction $\widetilde{v}_{k}^{(k)}$ is obtained from $v_{k}^{(k)}$ 
by removing one box from the part $E_k^{(k)}$. By definition of $\Rlh^{-1}$ and $\Rlb^{-1}$, 
$\widetilde{v}_{k+1}^{(k)}=\widetilde{v}_{k}^{(k)}$ and $v_{k+1}^{(k)}=v_{k}^{(k)}$. Moreover, 
we have $\widetilde{v}_{k+1}^{(j)}=\widetilde{v}_{k}^{(j)}$ and $v_{k+1}^{(j)}=v_{k}^{(j)}$ for
all $j \le k$. Therefore, the difference between $cJ_{k+1}^{(k)}$ and $\widetilde{cJ}_{k+1}^{(k)}$
is determined by the difference between $cJ_k^{(k)}$ and $\widetilde{cJ}_k^{(k)}$, the change of
sequence of rectangles by $\Rlb^{-1}\circ\Rlh^{-1}$ and the location of the added box 
$E_{k+1}^{(k+1)} \setminus D_{k+1}^{(k+1)}$. We 
note that the sequence of rectangles change decreases the contribution to $cJ_{k+1}^{(k)}(s)$ by $1$
for any string $s$, while the added box $E_{k+1}^{(k+1)} \setminus D_{k+1}^{(k+1)}$ increases the 
contribution to $cJ_{k+1}^{(k)}(s)$ for all $s$ of size $>|D_{k+1}^{(k+1)}|$. Thus 
\begin{itemize} 
   \item for strings $s \le D_{k+1}^{(k+1)}$, $\widetilde{cJ}_{k+1}^{(k)}(s)=cJ_{k+1}^{(k)}(s)$;
   \item for strings $D_{k+1}^{(k+1)} < s \le D_{k}^{(k)}$, 
         $\widetilde{cJ}_{k+1}^{(k)}(s)=cJ_{k+1}^{(k)}(s)-1$;
   \item for strings $D_{k}^{(k)} < s \le D_{k-1}^{(k-1)}$, 
         $\widetilde{cJ}_{k+1}^{(k)}(s)=cJ_{k+1}^{(k)}(s)+1$;
   \item for strings $s >D_{k-1}^{(k-1)}$, $\widetilde{cJ}_{k+1}^{(k)}(s)=cJ_{k+1}^{(k)}(s)$. 
\end{itemize} 

By the second bullet point above $E_k^{(k)}$ is the smallest singular string in 
$(v_{k+1},cJ_{k+1})^{(k)}$ is greater or equal to $E_{k+1}^{(k+1)}$. 

Consider the difference between $(v_{k+1},cJ_{k+1})^{(j)}$ and 
$(\widetilde{v}_{k+1},\widetilde{cJ}_{k+1})^{(j)}$ for $1 \le j < k-1$. 
Notice that $(v_{k+1},cJ_{k+1})^{(j)}=(v_k,cJ_k)^{(j)}$ and 
$(\widetilde{v}_{k+1},\widetilde{cJ}_{k+1})^{(j)}=(\widetilde{v}_k,\widetilde{cJ}_k)^{(j)}$.
Thus by induction we have: $\widetilde{v}_{k+1}^{(j)}$ removes one box from the part 
$E_j^{(j)}$ in $v_{k+1}^{(j)}$, and $\widetilde{cJ}_{k+1}^{(j)}$ is such that
\begin{itemize} 
   \item the shortened part has colabel 0;
   \item for strings $s\le |D_{j+1}^{(j+1)}|$, $\widetilde{cJ}_{k+1}^{(j)}(s)=cJ_{k+1}^{(j)}(s)$;
   \item for strings $D_{j+1}^{(j+1)}<s \le D_{j}^{(j)}$, 
         $\widetilde{cJ}_{k+1}^{(j)}(s)=cJ_{k+1}^{(j)}(s)-1$;
   \item for strings $D_{j}^{(j)}<s \le |D_{j-1}^{(j-1)}|$, 
         $\widetilde{cJ}_{k+1}^{(j)}(s)=cJ_{k+1}^{(j)}(s)+1$;
   \item for strings $s>D_{j-1}^{(j-1)}$, $\widetilde{cJ}_{k+1}^{(j)}(s)=cJ_{k+1}^{(j)}(s)$. 
\end{itemize} 

By the second bullet point above $E_j^{(j)}$ is the smallest singular string in 
$(v_{k+1},cJ_{k+1})^{(j)} \ge E_{j+1}^{(j+1)}$. 
\end{proof}

\begin{lemma}
The diagram \eqref{eq: triangle} commutes.
\end{lemma}
\begin{proof}
By Lemma~\ref{lemma: alpha}, the difference between $\Phi(S_{t-1}\otimes q)=(v_t,cJ_t)$ and 
$\Phi((\Pslide(S))_t\otimes q)= \Phi(\Pslide(S \otimes q))=(\widetilde{v}_t,\widetilde{cJ}_t)$ is  
\begin{enumerate}  
  \item for $j > t$, $(\widetilde{v},\widetilde{cJ})^{(j)}=(v,cJ)^{(j)}$. 
  \item for $j = t$, $\widetilde{v}^{(t)}=v^{(t)}$, $\widetilde{cJ}^{(t)}$ is obtained from $cJ^{(t)}$ by
        \begin{itemize} 
           \item for strings $s \le D_{t-1}^{(t-1)}$, $\widetilde{cJ}^{(t)}(s)=cJ^{(t)}(s)+1$;
           \item for strings $s> D_{t-1}^{(t-1)}$, $\widetilde{cJ}^{(t)}(s)=cJ^{(t)}(s)$. 
        \end{itemize} 
  \item for $j = t-1$, $\widetilde{v}^{(t-1)}$ removes one box from the part 
        $E_{t-1}^{(t-1)}$ in $v^{(t-1)}$. $\widetilde{cJ}^{(t-1)}$ is such that
        \begin{itemize} 
           \item the shortened part has colabel 0;
           \item for strings $s\le D_{t-1}^{(t-1)}$, $\widetilde{cJ}^{(t-1)}(s)=cJ^{(t-1)}(s)-1$;
           \item for strings $D_{t-1}^{(t-1)}< s \le D_{t-2}^{(t-2)}$, 
                 $\widetilde{cJ}^{(t-1)}(s)=cJ^{(t-1)}(s)+1$; 
           \item for strings $s>D_{t-2}^{(t-2)}$, $\widetilde{cJ}^{(t-1)}(s)=cJ^{(t-1)}(s)$. 
        \end{itemize} 
  \item for $1 \le j < t-1$, $\widetilde{v}^{(j)}$ removes one box from the part 
        $E_j^{(j)}$ in $v^{(j)}$. $\widetilde{cJ}^{(j)}$ is such that
        \begin{itemize} 
           \item the shortened part has colabel 0;
           \item for strings $s\le D_{j+1}^{(j+1)}$, $\widetilde{cJ}^{(j)}(s)=cJ^{(j)}(s)$;
           \item for strings $D_{j+1}^{(j+1)}<s \le D_{j}^{(j)}$, 
                 $\widetilde{cJ}^{(j)}(s)=cJ^{(j)}(s)-1$;
           \item for strings $D_{j}^{(j)}<s \le D_{j-1}^{(j-1)}$, 
                 $\widetilde{cJ}^{(j)}(s)=cJ^{(j)}(s)+1$;
           \item for strings $s>D_{j-1}^{(j-1)}$, $\widetilde{cJ}^{(j)}(s)=cJ^{(j)}(s)$. 
        \end{itemize} 
\end{enumerate}  

Let $(u,cI)=\Phi(S\otimes q)=\Rlb^{-1}\circ\Rlh^{-1}(\Phi(S_{t-1} \otimes q),n+2)$. 
Thus, $u$ can be obtained from $v$ by adding a box to $D_{t}^{(j)}$ for 
$j=n+1$ down to $t$. Now we use the fact that $n+2$ does not appear anywhere in $q$, 
which implies that $v^{(n+1)}$ is empty. Hence $D_t$ is the sequence of empty strings. This implies 
that the colabels of all unchanged strings are preserved when passing from 
$(v,cJ)$ to $(u,cI)$ (since the vacancy number for all unchanged strings are preserved,
and their labels are preserved). Then the difference between $u$ and 
$\widetilde{v}$ is given by the following:
\begin{enumerate}  
  \item for $j > t$, $(\widetilde{v},\widetilde{J})^{(j)}$ is obtained from $u$ by removing
        a string of size 1; 

  \item for $j = t$, $\widetilde{v}^{(t)}$ is obtained from $u^{(t)}$ by removing a box from
        the part $E^{(t)}_t$; 

  \item for $j = t-1$, $\widetilde{v}^{(t-1)}$ is obtained from $u^{(t-1)}$ by removing a box
        from the part $E^{(t-1)}_{t-1}$; 
         
  \item for $1 \le j < t-1$, $\widetilde{v}^{(j)}$ is obtained from $u^{(t)}$ by removing a box
        from the part $E^{(j)}_{j}$; 
\end{enumerate}  
By Lemma~\ref{lemma: alpha}, the sequence of boxes removed is precisely the 
$\Rslide$-sequence. Furthermore, the difference between $(v,cJ)$ and 
$(\widetilde{v},\widetilde{cJ})$ mandated by $\alpha$ precisely makes 
$\Rslide(u,cI)=(\widetilde{v},\widetilde{cJ})$.
\end{proof}

\appendix
%%%%%%%%%%%%%%%%%%%%%%%%%%%%%%%%%%%%%%%%%%%%%%%%%%%%%%%%%%%%%%%%%%%%%%%%%%%%%%%%%%%%%%%%%%%%%%%%%%%%%%%%%%%%%%%%%%%%%%%%%
%APPENDIX 
%%%%%%%%%%%%%%%%%%%%%%%%%%%%%%%%%%%%%%%%%%%%%%%%%%%%%%%%%%%%%%%%%%%%%%%%%%%%%%%%%%%%%%%%%%%%%%%%%%%%%%%%%%%%%%%%%%%%%%%%%
\section{Proof of Theorem~\ref{thm: Psi=Phi}} \label{pf: Psi=Phi}

The aim here is to prove Theorem~\ref{thm: Psi=Phi}. To do this, we will actually prove a
stronger statement. Let us first generalize the construction of $\Psi$ in 
Section~\ref{pf: base case 1}.

Let $c \ge 0$, and $0< t \le r$. Let $p \in \Path_n((t,1),(r,c))$ and write
\begin{equation*} 
   \Yboxdim25pt
p = T \otimes S = 
\young(\Toneone,\vdots,\Ttone)
       \otimes
\young(\Loneone\cdots\cdots\Lonec,%
       \vdots\vdots\vdots\vdots,%
       \Lrone\Lrtwo\cdots\Lrc)  \; .
\end{equation*}
We require that $T_{k,1} \le S_{k,1}$ for $k = 1,\ldots, t$, and write $p$ as
\begin{equation*} 
   \Yboxdim25pt
p = 
\young(\Toneone\Loneone\cdots\cdots\Lonec,%
       \vdots\vdots\vdots\vdots\vdots,%
       \Ttone\vdots\vdots\vdots\vdots,%
       :\vdots\vdots\vdots\vdots,%  
       :\Lrone\Lrtwo\cdots\Lrc)  \; .
\end{equation*} 
Define $\psi(p)=(u_1,\ldots,u_n)$ where each $u_k$ is the area of $p$ that is formed by all boxes that
\begin{enumerate} 
  \item are on the $j$-th row for $j \le k$;
  \item contain numbers $>k$. 
\end{enumerate} 

\begin{example}
Let $p \in \Path_7((2,1),(4,4))$ be $\young(1234,2445,:666,:778)$, then 
\begin{equation*}
\psi(p)\;=\;
(\young(234),\young(:34,445),\young(::4,445,666),\young(::5,666,778),\young(666,778),\young(778),\young(8))\; .
\end{equation*}
\end{example}
Remark~\ref{rmk: area} also applies here.

\begin{definition}
Define $\Psit(p)=(\psi(p),J) \in \RC_n((t,1),(r,c))$, where $\psi(p)$ is viewed as a sequence of partitions,
and $J$ is such that it sets the colabels for all parts in $u_t$ of size $\le c$ to 1 and 
sets the colabels of all other parts to 0.
\end{definition}

\begin{lemma}
For $p$ given as above, $\Phi(p)=\Psit(p)$.
\end{lemma}
\begin{proof}
We prove this statement by induction.

For a fixed $r \ge 0$, let $E_r(c,t)$ be the following statement with $c\ge 0$ and $0\le t\le r$ as 
free variables:
"For any $p = T \otimes S \in \Path_n((t,1),(r,c))$ with $T_{k,1} \le S_{k,1}$ for
$k=1,\ldots,t$ we have $\Phi(p)=\Psit(p)$."

The induction is on $E_r(c,t)$ with $(c,t)$ in the lattice $\mathbb{Z}_{\ge 0}\times [r]$.
For the base case, $E_r(0,0)$ is true since both $\Phi(p)$ and $\Psit(p)$ are lists of empty partitions.
The induction step has following two cases:
\begin{enumerate}
  \item Assume $E_r(c,t)$ for $c \ge 0$ and $t < r$, show $E_r(c,t+1)$;
  \item Assume $E_r(c,r)$ for $c \ge 0$, show $E_r(c+1,0)$. 
\end{enumerate}
For the first case, given $\Phi(p)=\Psit(p)$ for the following 
\begin{equation*} 
   \Yboxdim25pt
p = 
\young(\Toneone\Loneone\cdots\cdots\Lonec,%
       \vdots\vdots\vdots\vdots\vdots,%
       \Ttone\vdots\vdots\vdots\vdots,%
       :\vdots\vdots\vdots\vdots,%  
       :\Lrone\Lrtwo\cdots\Lrc)  
\end{equation*} 
we would like to compare $\Phi(p')$ and $\Psi(p')$ for 
\begin{equation*} 
   \Yboxdim25pt
p' = 
\young(\Toneone\Loneone\cdots\cdots\Lonec,%
       \vdots\vdots\vdots\vdots\vdots,%
       \Ttone\vdots\vdots\vdots\vdots,%
       a\vdots\vdots\vdots\vdots,%  
       :\Lrone\Lrtwo\cdots\Lrc)  
\end{equation*} 
where $S_{t+1,1} \ge a > \Ttone$.

The change from $\Phi(p)\in\RC((t,1),(r,c))$ to $\Phi(p')\in\RC((t+1,1),(r,c))$ 
is caused by $\Rlb^{-1}\circ\Rlh^{-1}$ which is described 
in the following algorithm: Let $s^{(a)}=\infty$. For $k=a-1$ down to $t+1$, select the 
longest singular string in $\Phi(p)^{(k)}$ of length $s^{(k)}$ (possibly of zero length) such that 
$s^{(k)} \le s^{(k+1)}$. Add a box to each of the selected strings, and reset their labels to make them 
singular with respect to the new vacancy number, leaving all other strings unchanged. 

By the construction of $\Psit(p)$, the inductive assumption that $\Phi(p)=\Psit(p)$, and 
$S_{t+1,1} \ge a$, we can conclude that $s^{(k)}=c$ for $k=a-1$ down to $t+1$ in the construction
of $\Phi(p')$. The changing of sequence of rectangles from $((t,1),(r,c))$ to $((t+1,1),(r,c))$ 
causes the colabels of all parts in $\Phi(p')^{(t+1)}$ of size $\le c$ set to 1 
(increased by 1 from 0), and the colabels of all parts in $\Phi(p')^{(t)}$ of size $\le c$ 
set to 0 (decreased by 1 from 1). This is precisely the effect of going from $\Psit(p)$ to 
$\Psit(p')$.
 
For the second case, going from $\Phi(p)$ to $\Phi(p')$ has the effect of $\Rls^{-1}$, which decreases
the colabels for all parts in $\Phi(p)^{(r)}$ of size $\le c$ by 1. Again this is precisely the result
of going from $\Psit(p)$ to $\Psit(p')$.
\end{proof}

\section{Several useful facts} \label{facts}
In this section, several facts that are repeatedly used in Section~\ref{pf: Rslide and Rls}
are stated and proved.

For any $p \in \Path_n$, let $\rc = \Phi(p)$. Then to each numbered box in $p$ one can
associate a $\Rlh^{-1}$-sequence.
In the case 
\begin{equation*} 
   \Yboxdim25pt
p= S \otimes T \otimes q = 
\young(\bONE,\vdots,\bJ)
       \otimes
\young(\Tonek\cdots\Toneone,%
       \vdots\vdots\vdots,%
       \Trk\cdots\Trone)  \; 
       \otimes
q \in \Path_n((j,1),(r,s),B'),
\end{equation*}
we adopt the definitions of $\rc_{j,s}$ and $D_{j,s}$ as given in Section~\ref{pf: Rslide and Rls}.
Unlike those in Section~\ref{pf: Rslide and Rls}, however, results in this section are general 
facts about rigged configurations and not just about $\Phi(\Dom(\Pslide))$.  

\begin{lemma} \label{lemma: same column}
For a fixed column $j$ of $T$ and for $1\le a<b \le r$, 
we have $D_{b,j}^{(k)} \le D_{a,j}^{(k)}$ for any $k \in [n]$.
\end{lemma}
\begin{proof}
By the definition of $\Rlh^{-1}$ (Definition~\ref{def: Rl_inverse}), we can inductively show 
$D_{a+1,j}^{(k)} \le D_{a,j}^{(k-1)}$ for all $k$ from $n$ down to $1$. Combining this with the fact that 
$D_{a,j}^{(k-1)} \le D_{a,j}^{(k)}$, we have $D_{a+1,j}^{(k)} \le D_{a,j}^{(k)}$. 
\end{proof}

\begin{lemma} \label{lemma: same row}
For a fixed row $i$ of $T$ and for $1\le c<d \le s$, 
we have $D_{i,d}^{(k)} > D_{i,c}^{(k)}$ for all $k$ such that $D_{i,c}^{(k)}\not = \infty$.
In the case $D_{i,c}^{(k)}= \infty$, $D_{i,d}^{(k)}$ must be also $\infty$.
\end{lemma}
\begin{proof}
We use the convention in this proof that $\infty$ plus any constant is $\infty$.
 
We proceed by induction on the row index $i$, on the following inductive hypothesis 
\begin{hypothesis}
For any $k \in [n]$, $\rc_{i-1,s+1}^{(k)}$ has a singular string of size 
$D_{i,s}^{(k)}+1$. Therefore $D_{i,s+1}^{(k)} \ge D_{i,s}^{(k)}+1$.  
\end{hypothesis}

By definition of $\Phi$, $\rc_{i-1,s+1}$ can be obtained from $\rc_{i,s}$ in 3 steps: 
\begin{enumerate}
   \item[S1.] From $\rc_{i,s}$ to $\rc_{r,s}$ by a 
             sequence of $\Rlb^{-1}\circ\Rlh^{-1}$ operations:
             \begin{equation*} 
                \rc_{r,s}= \Rlb^{-1}\circ\Rlh^{-1} 
                        (\cdots\Rlb^{-1}\circ\Rlh^{-1}
                           (\rc_{i,s},T_{i+1,s})\cdots,T_{r,s}).
             \end{equation*} 
   \item[S2.] From $\rc_{r,s}$ to $\rc_{0,s+1}$ by $\rc_{0,s+1}=\Rls^{-1}(\rc_{r,s})$.
   \item[S3.] (This step is empty for $i=1$.) From $\rc_{0,s+1}$ to $\rc_{i-1,s+1}$ by a 
             sequence of $\Rlb^{-1}\circ\Rlh^{-1}$ operations:
             \begin{equation*} 
                \rc_{i-1,s+1}= \Rlb^{-1}\circ\Rlh^{-1} 
                        (\cdots\Rlb^{-1}\circ\Rlh^{-1}
                           (\rc_{i,s},T_{1,s+1})\cdots,T_{i-1,s+1}).
             \end{equation*} 
\end{enumerate}

By definition of $\Rlh^{-1}$, for any $k \in [n]$, $\rc_{i,s}^{(k)}$ has a 
singular string of size $D_{i,s}^{(k)}+1$. By Lemma~\ref{lemma: same column}, all the 
consecutive applications of $\Rlb^{-1}\circ\Rlh^{-1}$ in S1 do not affect 
the singularity of the above strings. Thus for any $k \in [n]$, $\rc_{r,s}^{(k)}$ has a 
singular string of size $D_{i,s}^{(k)}+1$. In S2, $\Rls^{-1}$ never makes any 
singular string into non-singular. Thus for any $k \in [n]$, $\rc_{0,s+1}^{(k)}$ has a 
singular string of size $D_{i,s}^{(k)}+1$. Finally, by induction and 
Lemma~\ref{lemma: same column}, we see that $D_{k,s+1}$ for any $k<i$ will not affect 
the singularity of the above strings. Thus in $\rc_{i-1,s+1}^{(k)}$ for each $k$, 
there is a singular string of size $D_{i,s}^{(k)}+1$. 

Note that $T_{i,s+1} \le T_{i,s}$. This implies
$D_{j,s+1} > D_{j,s}$ for all $j \ge T_{i,s+1}$. Using this as base case, and the result
we just proved above, a downward induction shows that $D_{i,s+1}^{(k)} \ge D_{i,s}^{(k)}$ 
for all $k\in [n]$ and as long as $D_{i,s}^{(k)}<\infty$, $D_{i,s+1}^{(k)} > D_{i,s}^{(k)}$. 
\end{proof}


\begin{thebibliography}{99} 

\bibitem{DS:2006}
L.~Deka, A.~Schilling,
\textit{New fermionic formula for unrestricted Kostka polynomials},
J. Combinatorial Theory, Series A  \textbf{113} (2006) 1435--1461.

\bibitem{Young Tableaux}
William Fulton
\textit{Young Tableaux}
London Mathematical Society Student Texts 35.

\bibitem{HK:2002}
J.~Hong, S.-J.~Kang,
Introduction to quantum groups and crystal bases,
Graduate Studies in Mathematics, \textbf{42}, American Mathematical Society, Providence, RI, 2002. xviii+307 pp. 

\bibitem{HT:2003}
F.~Hivert, N.~M.~Thi\'ery,
\textit{MuPAD-Combinat, an Open-Source Package for Research in Algebraic Combinatorics},
S\'eminaire Lotharingien de Combinatoire \textbf{51} (2003) [B51z] (70 pp).\newline
{\tt http://mupad-combinat.sourceforge.net/}

\bibitem{KKR:1986}
S.~V.~Kerov, A.~N.~Kirillov and N.~Yu.~Reshetikhin
\textit{Combinatorics, the Bethe ansatz and representations of the symmetric group}, 
Zap.Nauchn. Sem. (LOMI) {\bf 155} (1986) 50--64.
(English translation: J. Sov. Math. {\bf 41} (1988) 916--924.)

\bibitem{KR:1988} 
A.~N.~Kirillov and N.~Y.~Reshetikhin,
\textit{The Bethe Ansatz and the combinatorics of Young tableaux},
J. Soviet Math. \textbf{41} (1988) 925--955.

\bibitem{KSS:2002}
A.N.~Kirillov, A.~Schilling, M.~Shimozono,
\textit{A bijection between Littlewood-Richardson tableaux and
rigged configurations},
Selecta Mathematica (N.S.) \textbf{8} (2002) 67--135.

\bibitem{KS:2009}
A.~Kuniba, R.~Sakamoto,
\textit{Combinatorial Bethe ansatz and generalized periodic box-ball system},
Rev. Math. Phys.  \textbf{20}  (2008),  no. 5, 493--527.

\bibitem{KOSTY:2006}
A.~Kuniba, M.~Okado, R.~Sakamoto, T.~Takagi, Y.~Yamada,
\textit{Crystal interpretation of Kerov-Kirillov-Reshetikhin bijection},
Nucl. Phys. B \textbf{740} (2006) 299--327.

\bibitem{Sage}
Sage, Open Source Mathematics Software, {\tt http://www.sagemath.org/}.

\bibitem{Sa:2008}
R.~Sakamoto,
\textit{Crystal interpretation of Kerov-Kirillov-Reshetikhin bijection II. Proof for $sl_n$ Case},
J. Algebr. Comb. \textbf{27} (2008) 55--98.

\bibitem{Sa:2009}
R.~Sakamoto,
\textit{Kirillov--Schilling--Shimozono bijection as energy functions of crystals},
International Mathematics Research Notices (2009) 2009, no. 4, 579-614.

\bibitem{Sa:2009a}
R.~Sakamoto,
\textit{Finding rigged configurations from paths},
RIMS Kokyuroku Bessatsu \textbf{B11} (2009) 1--17.

\bibitem{S:2006}
A.~Schilling,
\textit{Crystal structure on rigged configurations},
International Mathematics Research Notices, Volume 2006, Article ID 97376, Pages 1--27.

\bibitem{Sch:1972}
M.~P.~Sch{\"u}tzenberger,
\textit{Promotion des morphismes d'ensembles ordonn\'es},
Discrete Math. \textbf{2} (1972) 73--94

\bibitem{SW:1999}
A.~Schilling, S.O.~Warnaar,
\textit{Inhomogeneous lattice paths, generalized Kostka polynomials and A$_{n-1}$ 
supernomials},
Commun. Math. Phys. \textbf{202} (1999) 359--401. 

\bibitem{Sh:2002}
M.~Shimozono,
\textit{Affine type $A$ crystal structure on tensor products of rectangles, Demazure 
characters, and nilpotent varieties},
J. Algebraic Combin. \textbf{15} (2002), no. 2, 151--187.

\bibitem{WQ:FPSAC2009}
Q.~Wang,
\textit{A promotion operator on rigged configurations}, 
Discrete Math. and Theor. Computer Sci. proc. AK (2009) 887--898.

\end{thebibliography}
\end{document}